\newcommand{\C}{\mathbb C}
\newcommand{\N}{\mathbb N}
\newcommand{\R}{\mathbb R}
\newcommand{\Z}{\mathbb Z}
\newcommand{\Id}{\mathrm{Id\,}}

\newcommand{\ga}{\gamma}
\newcommand{\si}{\sigma}
\newcommand{\la}{\lambda}
\newcommand{\Ker}{\mathrm{Ker\,}}

\newcommand{\Ind}{\mathrm{Ind\,}}
\newcommand{\J}{J\,}
\def\f{\Phi(X,Y)}
\def\fz{\Phi_{0}(X,Y)}

\def\01{[0,1]}
\def\!{^{-1}}

\newcommand{\cd}{\mathcal{D}}
\newcommand{\cf}{\mathcal{F}}
\newcommand{\co}{\mathcal{O}}

\newcommand{\cg}{\mathcal{G}}
\newcommand{\cL}{\mathcal{L}}
\newcommand{\cLL}{\mathcal{L_\lambda}}

\newcommand{\cB}{\mathcal{B}}

\newcommand{\cq}{\mathcal{Q}} 
\newcommand{\cS}{\mathcal{S}}
\newcommand{\cSL}{\mathcal{S}_\la}
\newcommand{\cP}{\mathcal{P}}
\newcommand{\cH}{\mathcal{H}}
\newcommand{\hkr}{H^{k+s}(\Omega;\R^m )}
\newcommand{\hk}{H^{k+s}(\Omega;\C^m )}
\newcommand{\hs}{H^{s}(\Omega;\C^m )}
\newcommand{\hsr}{H^{s}(\Omega;\R^m )}
\newcommand{\hsw}{H^{s}_W(\co;\C^m )}

\newcommand{\hck}{H_{comp}^{k+s}(\co;\C^m )}

\newcommand{\hl}{H_{loc}^{s}(\co;\C^m )}

\newcommand{\hsb}{H^{s}(\Omega;\C^m)\times H^+(\partial\Omega;\C^r )}
\newcommand{\hsbr}{H^{s}(\Omega;\R^m)\times H^+(\partial\Omega;\R^r )}
\newcommand{\hbb}{H^+(\partial \Omega;\C^r )}
\newcommand{\hbbr}{H^+(\partial \Omega;\R^r )}
\newcommand{\hkn}{H^{k+s}(\R^n;\C^m)}

\newcommand{\hsn}{H^{s}(\R^n;\C^m )}

\newcommand{\ind}{\operatorname{ind}}
\newcommand{\sgn}{\operatorname{sgn}}
\newcommand{\im}{\operatorname{Im}}

\newcommand{\coker}{\operatorname{coker}}
\newcommand{\hin}{\hbox{\,in\, }}
\newcommand{\hfor}{\hbox{ \,for\, }}
\newcommand{\hand}{\hbox{ \,and\, }}
\newcommand{\hif}{\hbox{ \,if\,} }
\newcommand{\hof}{\hbox{\,of\,} }

\newcommand{\hforall}{\hbox{ \,for\  all\, }}
\newcommand{\LL}{\Lambda}

\newcommand{\pc}{\,Ell(\R^n)\,}
\newcommand{\ch}{\operatorname{ch}}
\newcommand{\he}{\operatorname{\mathbb H}^{ev}}
\newcommand{\h}{\operatorname{\mathbb H}}
\newcommand{\bt}{\begin{theorem}}
\newcommand{\et}{\end{theorem}}
\newcommand{\bc}{\begin{corollary}}
\newcommand{\ec}{\end{corollary}}
\newcommand{\bp}{\begin{proposition}}
\newcommand{\ep}{\end{proposition}}
\newcommand{\bl}{\begin{lemma}}
\newcommand{\el}{\end{lemma}}
\newcommand{\br}{\begin{remark}}
\newcommand{\er}{\end{remark}}
\newcommand{\bd}{\begin{definition}}
\newcommand{\ed}{\end{definition}}
\newcommand{\be}{\begin{equation}}
\newcommand{\ee}{\end{equation}}
\newcommand{\ra}{\rightarrow}
\newcommand{\eqr}{\eqref}
\documentclass[thmsa,a4paper]{amsart}

\usepackage{amsmath}
\usepackage{amsfonts}
\usepackage{amssymb}
\input{diagxy}

\begin{document}
\newtheorem{theorem}{Theorem}[section]
\newtheorem{corollary}[theorem]{Corollary}
\newtheorem{lemma}[theorem]{Lemma}
\newtheorem{proposition}[theorem]{Proposition}
\newtheorem{remark}{Remark}[section]
\newtheorem{definition}{Definition}[section]
\numberwithin{equation}{section}
\numberwithin{theorem}{subsection}
\numberwithin{definition}{subsection}
\numberwithin{remark}{subsection}
\begin{abstract}We associate to a parametrized family $f$ of nonlinear Fredholm maps possessing  a trivial branch of zeroes  an {\it index of bifurcation} $\beta(f)$  which provides an algebraic measure for the number of bifurcation points from the trivial branch. The index $\beta(f)$ is derived from the index bundle of the linearization of the family along the trivial branch by means of the generalized $J$-homomorphism.  Using the Agranovich reduction and a cohomological form of the Atiyah-Singer family index theorem, due to Fedosov, we compute the bifurcation index of a multiparameter  family of nonlinear elliptic boundary value problems from the principal symbol of the linearization along the trivial branch. In this way we obtain criteria for bifurcation of solutions of nonlinear elliptic equations which cannot be achieved using the classical Lyapunov-Schmidt method.
\end{abstract}
\title{Bifurcation of Fredholm maps I;
The Index Bundle and Bifurcation}
\author{J. Pejsachowicz }

\address{Dipartimento di Matematica Politecnico di Torino,\\
    Corso Duca degli Abruzzi  24, 
       10129 Torino, Italy.}
\email{jacobo.pejsachowicz@polito.it}



\date{\today}

\thanks{This work was supported by  MIUR-PRIN2007-Metodi variazionali e topologici nei fenomeni nonlineari}
\maketitle
\tableofcontents 
\date{16/10/2009}
\specialsection{\footnotesize INTRODUCTION AND STATEMENTS OF  THE MAIN RESULTS }\label{sec:1}

\vskip10pt

\subsection{Introduction}\label{secá0}
The main purpose of the article is to present a comprehensive  account of  the relationship between elliptic topology  and bifurcation theory. More precisely,  between the index bundle of a family of linear  Fredholm operators  and  bifurcation of solutions of nonlinear elliptic equations from a trivial branch. 

 Bifurcation from a trivial branch is one of the oldest notions  of bifurcation in mathematics.  Roughly speaking, the scheme is as follows: assuming  that there is a known (trivial) branch of solutions of a parametrized family of problems,  find necessary and sufficient conditions for the appearance of nontrivial solutions arbitrary close to some points (called bifurcation points) of the trivial branch.  The above framework arises in several fields  belonging to pure and applied mathematics, which explains the interest in the formulation of a structured theory going  beyond a collection of examples.  
 
 Although the first studies of specific bifurcation phenomena can be traced back  to Euler and Jacobi,  bifurcation theory  was born with Poincar\'e as a special  chapter of  his qualitative  theory of dynamical systems. The most important tool for the analysis of bifurcation from a trivial branch is  the  Lyapunov-Schmidt reduction, which  leads a given bifurcation problem for  integral and differential equations to a locally equivalent problem for a finite number of nonlinear equations in a finite number of indeterminates. 
 
 Bifurcation can arise only at singular points of the linearization at the trivial branch, i.e., points  belonging to the trivial branch at which the linearized operator in the normal direction to the branch fails to be invertible.  One of the typical assumptions of the Lyapunov-Schmidt method is that  singular points are isolated.   Assuming this, there is a large variety of methods which, combined with the  Lyapunov-Schmidt reduction, provide  criteria for the appearance of nontrivial solutions close to the singular point 
 \cite{[Bu-To], [Ki], [Sat], [Da-Sa-Ts],[Iz-Vi]}. 
 
  The choice of the approach  depends on the nature of the problem at hand. However, the most popular ones use either the singularity theory or  topological methods.   In the first case, whether the point under consideration is a bifurcation point or not, is solved investigating higher order jets of the reduced map.  In the topological approach, particularly useful in the several parameter case, the presence of bifurcation is determined from  topological  invariants, described in section \ref{sec:3.3}.  The books \cite{[Va-Tr],[Ch-Ha], [Go],[Ra]} are only few of the several possible  references to the first method.  Ize's  PhD thesis \cite{[Iz]}, his review \cite{[Iz-1]} together with \cite{[Al],[Ba-Kr-St]}  provide a good introduction to the second one.  
 
  In this paper we will consider bifurcation of parametrized families of Fredholm maps from a topological  viewpoint  which is different from the well established method mentioned above.   We will not make any assumption about the nature of singular points of the linearization but we will heavily rely on the nontrivial topology of the parameter space. More precisely,  we will  look for homotopy invariants of the family of linearizations at points of the trivial branch whose non-vanishing entails the presence of at least one bifurcation point. 
 
 It should be noted that  invariants of this type exist because the homotopy groups of the space of  linear Fredholm operators between infinite dimensional Banach spaces are nontrivial.  Thus, our theory is strongly tied to homotopy  theory of Fredholm operators, i.e., elliptic topology. On the other hand, it complements the local point of view developed by Alexander and Ize providing criteria for bifurcation that are different  from the ones that can be obtained using the Lyapunov-Schmidt reduction. 
 
  To some extent,  our approach  was inspired by the successful use of  elliptic invariants in handling various  linear PDE problems in geometry and analysis.  For example, in \cite{[Hi]} the index bundle for families  was used  with the purpose to find Riemannian manifolds such that  the dimension of the space of harmonic spinors varies with the metric. In \cite{[Va-Wh]} the same  method was applied  to  determine  spectral gaps  of Dirac operators. Several generalizations of  Lichnerowitz's  theorem relating the  $A$-genus  of a spin manifold to the non existence of a metric with positive scalar curvature are rooted on similar arguments.  Their basic idea is to evaluate the index bundle of the relevant family of linear Fredholm operators of index $0$ using family index theorems. If  the index bundle  is  nontrivial, then $\ker L_\la \neq \{0\},$ for  at least one value of the parameter $\la$.  What we will  show in this paper is that the above argument works for nonlinear Fredholm maps as well, but at the cost of introducing one extra tool:  the generalized $J$-homomorphism.  
 \vskip2pt
Our goals  are:
\vskip1pt

1) Given a family $\{f_\lambda\}_{\lambda \in \Lambda}$ of $C^{1}$-Fredholm maps depending continuously on a parameter belonging to a finite $CW$-complex $\Lambda$ such that $f_\la (0) =0 \hforall \lambda\in\Lambda,$  we will define an  \emph{index of bifurcation points}  $\beta(f)$  which, much in the same way as the Lefschetz number in fixed-point theory, provides an algebraic measure of the total number of bifurcation points of the family $f.$  The index $\beta(f) $ takes values in a finite group $J(\Lambda).$   It only depends on the homotopy class of the family $  \{ L_{\la } = Df_{\la}(0):  \la \in \Lambda\} $ of linearizations of $f$ at points of the trivial branch.  In particular, when  $f$ is defined by a family of nonlinear elliptic differential operators,  $\beta(f)$ depends only on the coefficients of leading terms of the linearization.  
\vskip1pt
2) We will introduce a local index of bifurcation $\beta(f,U),$  analogous to the local fixed-point index, which interpolates between $\beta(f)$ and the index at an isolated point  derived from the Alexander-Ize bifurcation invariant. It is defined only if $L_\la $ is invertible for $\la$ outside of a compact subset of $U$ and preserved by homotopies of this type. In the case of nonlinear elliptic differential operators, in general, $\beta(f,U)$  depends on lower order terms  of the  linearized equations as well.
\vskip1pt
 3)  For particular families of nonlinear elliptic boundary value problems parametri\-zed by $\R^q$ we will compute the index of bifurcation  from the principal  symbol of the linearization along  the trivial branch using the  Agranovich reduction, Atiyah-Singer family index theorem and known results about the generalized $\J$-homomorphism. In this way we will obtain sufficient conditions for the existence of nontrivial  solutions bifurcating from the trivial branch for nonlinear elliptic problems with general  boundary conditions of Shapiro-Lopatinskij type.  Finally, using the local index, we will obtain conditions for the existence of multiple bifurcation points.
  \vskip 3pt
 
 For families parametrized by $\R^q$ the results are particularly striking. While the proofs involve some amount  of algebraic topology,  the complete knowledge of the $J$-groups of spheres and Fedosov's  formula for the Chern character of the index bundle allows to state our main bifurcation result,  Theorem \ref{th:40},  in terms of divisibility of a number computed as an integral of a differential form constructed explicitly  from the  principal symbol of the linearization at the trivial branch.

 Let us remark  that due to the invariance  of $\beta(f)$  under lower order perturbations,  its  nonvanishing  provides stronger bifurcation results than the ones obtained using  the classical approach, which always need some knowledge of the solutions of the linearized equations. On the negative side one can say that, precisely for the same reasons, $\beta(f)$ frequently vanishes. For instance, when the leading coefficients of the linearization do not depend on the parameter.  In this case one has to resort  to the local index  in order to detect bifurcation points. Pushing the analogy with the  fixed-point theory one step further,  the role of the Atiyah-Singer formula in our theory is reminiscent of the role of the Lefschetz-Hopf formula there. 

\vskip1pt

In the case of semilinear Fredholm maps the proof of the main abstract result, Theorem \ref{th:1.1}  is  simpler, and was sketched in \cite{[Pe]}. Simple examples of a direct calculation of the bifurcation index from the data of the problem, using elementary family index theorems, can be found in \cite{[Fi-Pe-2]} and \cite{[Pe-3]}. The first deals with nonlinear Sturm-Liouville problems while the second studies bifurcation of homoclinic orbits.  

Here for the first time we deal with general nonlinear Fredholm maps and use the Atiyah-Singer theorem in order to compute the bifurcation index of a large family of elliptic boundary value problems with general boundary conditions.   Hence,  we will keep the presentation as complete and self-contained as possible. Taking into account the mixed nature of the subject, we will carefully introduce  the terminology used in the paper and prove most of the assertions. Some of our results  from chapters 2 and 3  were announced without proof in \cite{[Pe-2]}. 
\vskip1pt

 The paper is structured as follows:  precise statements of the results concerning  item $ i), \, 1\leq i\leq 3,$ of the above list are formulated in subsection $1.(i+1)$ of this section and proved together with some generalizations and corollaries in section  $i+1$  with the same title.  Subsection $1.5$  contains   several comments to related work and  eventual  further developments. There are three appendices.  In the first we sketch out  the proof of standard  properties of the index bundle. The second reviews some well  known results about Fredholm properties of  maps induced on Hardy-Sobolev spaces by linear and nonlinear elliptic operators.  The third is devoted to Fedosov's formula for the Chern character of the index bundle of a family of elliptic pseudo-differential operators. 
 
 Finaly,  I would like to thank  Ernesto Buzano,  Nils Waterstraat and  Victor Zviagin for their comments and  generous help.

 \vskip20pt
\subsection{Index bundle and the  index of bifurcation points}\label{sec:1.1}    

Let  $X,Y$  be real  Banach spaces,  $O$ be an open subset of $X,$  and let  $\Lambda$ be a  finite connected CW-complex.  A {\it family  of $C^n$-maps, $  0\leq n\leq \infty ,$  continuously parametrized by $\Lambda$ }  is a continuous map $f\colon \Lambda \times O\rightarrow X$  such that for each ${\la} \in  \Lambda$ the map  $f_{\la}\colon O\rightarrow X$ defined by $f_{\la}(x)= f(\la ,x) $ is  $C^n$   and,  for all $k\leq n,$ the $k$-th derivative of $f$ in direction $x,$   $D_x ^k f\colon \Lambda\times O \rightarrow L^k (X,Y),$ is  continuous in norm topology of the space of $L^k (X,Y)$ of $k$-forms on $X$ with values in $Y.$  
  
  Parametrized families of  $C^n$-maps are a particular case of fiberwise  $C^n$-maps, i.e.,  morphisms  in   the category  of  $C^n$-Banach manifolds over $ \Lambda$ (see \cite{[Cr], [At-Si-4]}). While  most of our arguments have a very  natural extension to this category, some problems arise  related to infinite dimensional structure groups.  Hence, we will consider  here only the product  case  $\Lambda \times O.$

We will deal mainly  with  families of $C^1$-Fredholm maps of index $0,$ which means that  $Df_{\la}(x)$ is a Fredholm operator of index $0$ for all $(\la ,x) \in  \Lambda \times O.$ We will further assume everywhere in this paper that  $O$ is  an open neighborhood of the origin  and  that $f(\la ,0)=0$  for all ${\la} $ in  $\Lambda .$ Solutions of the  equation $f(\la ,x)=0$ of the form $(\la ,0)$ are called \ {\it trivial.} \ The set $T= \Lambda \times \{0\}$  is called the {\it trivial branch.}  As  a rule we will identify  the trivial branch with  the parameter space  $\Lambda.$  
\begin{definition}\label{defbif}
A {\it bifurcation point  from the trivial branch} of solutions of  the equation $f(\la ,x)=0$ is a point ${\la} _*\in \Lambda$ such that every neighborhood of $(\la _*,0)$ contains nontrivial solutions of this equation. 
\end{definition}

In what follows, we will denote  with $\cL(X,Y)$ the  Banach space of all bounded operators from $X$ to $Y,$  with $\f$ (resp. $\Phi_k(X,Y)$ ) the open subspace of all Fredholm operators (resp.  those  of index $k).$ 

 The {\it  linearization of the family  $f$  along the trivial branch}  is the family of operators  $L\colon \Lambda\rightarrow \Phi_0(X,Y),$ where  $L_\la =Df_{\la}(0) $ is the Frechet derivative of  $f_{\la}$ at $0$. 
    
  Bifurcation  can only occur at  singular points of the linearization, i.e., the points ${\la} \in \Lambda$  such that  $\ker L_\la \neq 0.$  When $ \Lambda$ is a smooth manifold and $f$ is $C^1$ the necessity of this condition  follows immediately from the implicit function theorem.  It holds  in our slightly more general framework too. Indeed,  in a  small enough neighborhood of  a point $\nu$ such that  $L_\nu$ is  nonsingular,  the equation $f(\la ,x) =0$ is equivalent to $ x =L_{\la}^{-1} g(\la ,x)$  where $g(\la ,x)= f(\la ,x) - L_\la x.$  Since $g(\la ,x) =o(||x||),$ by the uniqueness of the fixed point of a  contraction the only solutions close to $(\nu,0)$ are the trivial ones. 

While necessary, the above condition is not sufficient for the appearance of nontrivial solutions close to the  given point of the trivial branch. Hence, in general, the set $Bif(f)$ of all bifurcation points of a family $f$ is only a  proper closed subset  of  the set  $ \Sigma(L)$  of all singular points of the  linearization $L$ along the trivial branch.  The purpose of the  {\it linearized bifurcation theory} is to obtain sufficient conditions for the existence of bifurcation points of $f$  in terms of the linearization $L.$  

Since bifurcation  arises only  at points of $\Sigma(L),$  the first topological  invariant that comes to mind is the obstruction to  deformation of  $L$ into a family without singular points. It is well known  that such an  obstruction is given by an  element of the reduced Grothendieck group of virtual vector bundles $\tilde{KO}(\Lambda)$, called {\it family index } or  {\it  index bundle}  \cite{[At],[Ja]} and denoted with $\Ind L$. However, since we are dealing  with nonlinear perturbations of  $L,$ we  have to take into account  the generalized $\J$-homomorphism  $\J\colon \tilde{KO}(\Lambda )\rightarrow \J(\Lambda )$ which associates to each  vector bundle  the  stable fiberwise homotopy class of  its unit sphere bundle. 

Quite naturally, our bifurcation invariant  is not $\Ind L$ but rather its image $\J(\Ind L)\in \J(\Lambda )$ under the generalized $J$-homomorphism.   In fact, we have:

\begin{theorem}
\label{th:1.1}
 Let  $f\colon \Lambda \times O\rightarrow Y$  be a family of $C^1$-Fredholm maps of index $0$ parametrized  by a connected finite $CW$-complex $\Lambda,$ such that  $f(\la,0)=0.$  If  $\Sigma(L)$ is a proper subset of $\Lambda $ and $ \beta (f) =\J(\Ind L)\neq 0,$  then  the family $f$ possesses at least one bifurcation point from the trivial branch.
  \end{theorem}
 
 The Stiefel-Whitney characteristic class $ \omega(E) = 1+\omega_1(E) + ...,  \omega_i(E) \in H^{i}(\Lambda;Z_{2})$ of a vector bundle $E$ over $\Lambda  $  is invariant under addition of trivial bundles  and hence  it is well defined on $\tilde{KO}(\Lambda).$  Moreover it factorizes through $J(\Lambda )$ because, by Thom's construction, it only depends on the stable fiberwise homotopy class of the  associated sphere bundle.   If $p$ is an odd prime, the same  holds for the total Wu class $ q(E) =  1+ q_1(E) + ..., q_{i}(E)  \in H^{2(p-1)i}(\Lambda;\Z_{p})$ \cite{[Mi]}. 
 In particular: 
\bc \label{stw}
Let $f$ and $\Sigma(L)$ be  as in the above theorem.  Then bifurcation arises if either $ \omega(\Ind L)\neq 1$ or  $ q(\Ind L) \neq 1$ for some odd prime $p.$\ec

The  nonvanishing of characteristic classes  of the index bundle of positive degree not only entails bifurcation but also gives some information about the size of the set $ Bif(f) $ of bifurcation points of $f$ and its position in the parameter space. We will study this in a companion paper \cite{[Pe-4]}.

\begin{remark}{\rm
  The assumption $\Sigma(L)\neq \Lambda$ can be relaxed (see section 2.4).  However, it is easy to see that nonvanishing of $\J(\Ind L)$ only,  does not  imply by itself  the existence of a bifurcation point. 
  
  For example, take a  family $L$ of Fredholm operators between Hilbert spaces whose kernels define a nonorientable bundle $\ker L$ over $\Lambda$ and such that  $\coker L$  is a trivial vector bundle. Families of ordinary differential operators with this property  can be found in \cite{[Pe-3]} and \cite{[Fi-Pe-2]}. By the above corollary,  $\J(\Ind L)\neq 0.$  Let  $Q$ and $Q'$ be projectors on   $\ker L$ and $\im L$ respectively, and let  $s$ be  a nowhere vanishing section of $ F =(\Id - Q')X \simeq  \coker L.$  Define the family  $f$  by  $f(\la, x) =L_\la x +||Q_\la x||^2 s (\la).$ Then the linearization of $f$ at the trivial branch is $L$ but $f$ has no bifurcation points.}
\end{remark}

 \vskip20pt
  
\subsection{A local index of bifurcation}\label{sec:1.2}
  
Let $U$  open subset of $ \Lambda,$  $f\colon U\times O \rightarrow Y$ be a family of $C^1$-Fredholm maps parametrized by $U$ such that $f(\la ,0)=0$ and let  $L$ be the  linearization of $f$  along the trivial branch.  A pair $(f,U)$ is called \emph{admissible} if the singular set $\Sigma(L) $ is a  compact,  proper  subset of $U.$  An admissible homotopy is a family of $C^1$-Fredholm maps parametrized by $[0,1]\times U$  such that  the set $$\Sigma(Dh)=\{ (t ,\la ) /  Dh_{(t ,\la )}(0)\,  \text{ is singular}\, \}$$  is  a compact subset of\  $ [0,1]\times U$ and $\Sigma(Dh_i); \, i=0,1$ are proper subsets of $U.$

Let us recall  that a Kuiper space is Banach space $Y$ such that the subspace $GL(Y)$ of all invertible operators in $\cL(Y)$  is  contractible.  

The main result in section \ref{sec:3}  is: 
 \medskip
\begin{theorem}
\label{th:locind}
Assume that  $Y$ is a Kuiper space. There exists a {\it local index of bifurcation}   which   assigns to each  admissible pair $(f,U)$ an element $\beta (f,U)\in \J(\Lambda )$
 verifying  the following properties:

\begin{itemize}
 	
 \item[$ \text{B}_1$]{ \rm Existence:} \  If $\beta(f,U)\neq 0,$ then the family $f$ has  a bifurcation point in $U$.

\item[$\text{B}_2$] {\rm Normalization:} \ 
 $\beta (f,\Lambda)=\beta(f) =\J(\Ind L)$. 

 \item[$\text{B}_3$]{ \rm Homotopy invariance:} \   If $h$ is an admissible homotopy, then 
 
  $$\beta (h_0,U)=\beta (h_1,U).$$

 \item[$\text{B}_4$] {\rm Additivity:} \  Let $(f,U)$ be admissible with $U\subset \bigcup  U_i$.
Put \ $\Sigma _i=\Sigma (f)\cap  U_i$ and $f_i= f\vert _{U_i}.$ \ \
 If  $\,\Sigma _i \cap \Sigma _j=\emptyset$\ and \ 
$\bigcup \Sigma _i=\Sigma (f),$ then $(f_i,U_i)$ are admissible and
 \[\beta (f,U)= \sum_{i}\beta (f_i,U_i).\] 
\item[$\text{B}_5$] { \rm Change of parameters:} \  Let
$\alpha \colon \Lambda' \rightarrow  \Lambda$ be a continuous map. Let $U$ be an open subset of $ \Lambda$ such that the pair  $(f,U)$ is admissible. If $U'=\alpha ^{-1}(U)$ and  
 $g\colon U'\times O\rightarrow Y$  is defined by $g(\la' ,x)=f(\alpha (\la'),x),$ then  $(g,U')$ is admissible and 
 $\beta (g,U')=\alpha^* \beta(f,U),$
where $\alpha^* \colon \J(\Lambda )\rightarrow \J(\Lambda')$ is the homomorphism induced by $\alpha$ in $\J$-groups. 
 
\item[$\text{B}_6$] {\rm  Isolated points:} \  Let $\la_0$ be an isolated point in $ \Sigma (L).$  Assume that  there exists 
 a neighborhood  $U$ of ${\la} _0$ homeomorphic to $\R^{n}$ 
such that $ \Sigma (L)\cap U = \{\la _0\}.$  Then,  identifying   $J(S^n)$  with image of the stable  $j$-homomorphism   $j \colon \pi_{n-1}GL(\infty)\rightarrow  \pi^s_{n-1},$ we have:  
 \begin{equation*}\label{eqn: jga}
 \beta(f,U)= q^* j(\gamma_f).
\end{equation*}
Here  $\gamma_f$ is the Alexander-Ize invariant (see section \ref{sec:3.3}),  $S^n$  is identified with  the one-point compactification  $U^+$ of $U,$  and  $q \colon \Lambda \rightarrow U^+ $ is the map collapsing  $\Lambda-U$ to the point at infinity.  
\end{itemize}
\end{theorem}
\begin{remark}
\label{rem:3}
{\rm  A special case of $\text{B}_4$ is the excision property: if 
 $(f,U)$
is admissible and $\Sigma \subset V\subset U,$ then $\beta (f,U)=\beta (f\vert_{V},V).$ 
It follows from this and $\text{B}_3)$ that $\beta (f,U)$ depends only on the germ 
of the family of linearizations  $L_{\la}=Df_{\la}(0)$ at $\Sigma$.}
\end{remark}

 Few words have to be said about the computation of  $J(\Lambda)$  since the bifurcation index takes values in this group. $J(S^q)$  has been completely determined in the seventies \cite{[Ad],[Hu]}. We will use this computation in the next subsection.  In order to state the result, let $\nu_p(s)$ denote the exponent to which the prime $p$ occurs in the prime decomposition of an integral number $s.$   Consider  the  number-theoretic function $m$ constructed as follows: the value  $m(s)$ is defined  through its prime decomposition by setting for  $p=2,$   $ \nu_2\left(m(s)\right)=2+\nu_2(s)$ if  $s\equiv 0 \mod 2$ and $ \nu_2\left(m(s)\right)=1$ if the opposite is true. While, if $p$ is an odd prime, then
  $\nu_p\left(m(s)\right)= 1 +\nu_p(s)$ if  $s\equiv 0 \mod (p-1)$  and $0$ in the remaining cases. In particular  $m(s)$ is always even.  With this said, 
  $\ J(S^q)=\Z_2 \hfor q \equiv 1\ \text {or}\ 2 \mod 8,$   $ J(S^q)=\Z_ {m(2s)}$ for $q=4s,$ and  is  trivial  in the remaining cases.
  
   The numbers $m(s)$  have a wide range  of distribution (see for example  \cite{[Al]}). However, what is important for us is that  the index of bifurcation $\beta(f,U)$ is  an integral $\mod m$  in the case $\Lambda=S^q.$   The same holds true for  $\Lambda =\mathbb{RP}^q,$  the real projective space.   For a finite CW-complex $\Lambda$  without two-torsion in homology the order of $J(\Lambda)$  can be estimated in terms of the homology of $\Lambda$ with coefficients in $J(S^q)$. 

 \vskip20pt

\subsection{Bifurcation of solutions of nonlinear elliptic BVP}\label{sec:1.3}
 
 In theorem \ref{th:40} below we  will state criteria for bifurcation of solutions of nonlinear elliptic boundary value problems in terms of the coefficients of the top order derivatives of  linearized equations.  In theorem   \ref{th:50} we will consider the existence of multiple  bifurcation points. 

 Let   $\Omega$ be  an open bounded subset of $\R^n$ with smooth boundary $\partial \Omega.$  Referring to the Appendix B for the notations, we will  consider  nonlinear 
 boundary value problems of the form 
 \begin{equation}\label{bvp1}
\left\{\begin{array}{l} \cf\,(\la, x,u,\ldots ,D^{k}u)=0 \hfor  x \in  \Omega , \\ \cg^i(\la, x,u,\ldots, D^{k_i} u)=0 \hfor  x \in  \partial \Omega , \, 1 \le  i \le r.\end{array}\right. \end{equation} 
Here,  $u\colon\bar\Omega\rightarrow \R^m$ is a vector function,  $\la\in \R^{q}$ is a parameter and, denoting with $k^*$ the number of $\alpha$'s  such that $ |\alpha| \leq k,$  $$ \cf \colon \R^q \times \bar \Omega \times  \R^{mk^*}  \rightarrow  \R^m \hand  \cg^{i} \colon \R^q \times \bar\Omega  \times  \R^{mk_i^*} \rightarrow  \R  $$  are smooth with  $\cf (\la,x, 0) = 0,\  \cg^i( \la,x,0)=0, \, 1 \le  i \le r .$
 
 We  also denote with $\cf$ the family of nonlinear differential operators  \[\cf\colon \R^{q}\times C^{\infty}(\bar\Omega; \R^{m}) \to C^{\infty}(\bar \Omega; \R^{m})\] induced by the map $\cf.$ 

The functions  $\cg^i$ define a family of nonlinear boundary operators  \[\cg\colon \R^{q}\times C^{\infty}(\bar\Omega; \R^{m}) \to C^{\infty}(\partial \Omega; \R^{r})\]  
\[ \cg(\la, x,u,\ldots ,D^{k}u)=(\tau \cg^1(\la, x,u,\ldots ,D^{k_{1}}u), \ldots, \tau \cg^r(\la, x,u,\ldots ,D^{k_{r}}u)),\] 
where $\tau$ is the restriction to the boundary.

We assume:

\begin{itemize}
 
 \item[$H_1)$]  For all $\la \in R^q,$  the linearization $(\mathcal L_\la(x,D), \mathcal B_\la(x,D))$  of $(\cf_\la, \cg_\la)$  at $u\equiv 0,$
is an elliptic boundary value problem in the sense of  definition \ref{EBVP} in Appendix B. 
\smallskip
 \item [$H_2)$]  The coefficients   $a^{ij}_\alpha, b^{ij}_\alpha $ of the linearization  $(\mathcal L, \mathcal B)$   extend  to  smooth functions  defined on  $S^q \times \bar \Omega,$ where  $S^q=\R^q\cup\{\infty\} $ is the one point compactification of $\R^q.$ 
    Moreover  the problem: 
  \[ \left\{ \begin{array}{l} {\mathcal L}_\infty (x,D) u(x)=\displaystyle \sum_{|\alpha|\leq k} a_{\alpha }(\infty, x) D^{\alpha }u(x) =f(x), \  x \in \Omega\\
\mathcal {B}^i_\infty (x,D) u(x)=\displaystyle \sum_{|\alpha|\leq k_i} b^i_\alpha (\infty,x)D^{\alpha }u(x)=g(x),\  x\in \partial\Omega, \, 1\leq i\leq r, \end{array} \right .\]  is elliptic and has a 
unique solution for every $f \in C^\infty(\bar \Omega; \R^m)$ and every  $g\in C^\infty (\partial \Omega; \R^r ).$
\smallskip 
   \item [$H_3)$]  \begin{itemize}\item[i)]  The coefficients $b^{ij }_\alpha(x), | \alpha|=k_i, \,1\leq i\leq r ,$ of the  leading terms of $ \mathcal B(\la, x,D)$ are independent of  $\la.$
        \item[ii)] There exist a compact set  $K\subset \Omega $  such that  the coefficients\\ $a^{ij}_{\alpha }(\la, x), \, |\alpha |=k, $ of the  leading terms   of   ${\mathcal L}_\la (x,D)$  are independent of $\la $ for $x\in \bar\Omega -K.$\end{itemize} 
\end{itemize}
 
Let us give a closer look to our assumptions. Since linear elliptic boundary value problems induce Fredholm operators on function spaces,  $H_1$  places the problem \eqref{bvp1}  in the framework of our abstract bifurcation theory applied to a family of nonlinear Fredholm maps $f.$ The assumption $H_2$  allows us  to compute the local bifurcation index $\beta(f,\R^n)$ from the index bundle of the extended family.  Finally,  $H_3$ is essential in order to carry out  the Agranovich reduction showing that  $\Ind L$  coincides with the index bundle of a family  $\cS$ of pseudo-differential operators whose principal symbol is the matrix function $\sigma$ defined in \eqref{symbo} below.                                

\vskip 5pt

 Let $p(\la,x,\xi ) \equiv \sum^{}_{\vert \alpha \vert =k}a_{\alpha }(\la,x)\xi ^{\alpha }$ be the principal symbol of $\cLL.$  Since  the symbol  is  defined in terms  of $D_j= -i \frac{\partial }{\partial x_j},$   $p(\la,x,\xi )$ is  a complex matrix  which verifies the {\it reality condition}  $p(\la,x,-\xi)=\bar p(\la,x,\xi).$  
 
  By ellipticity, $ p(\la,x,\xi )\in GL(m;\C)$ if $\xi\neq 0.$ On the other hand, by $H_3,$
    \[ p (\la,x,\xi )= p (\infty,x,\xi )  \hfor x \in \bar\Omega-K.\]
      Putting   \[\sigma(\la,x,\xi)= Id  \ \text{ for any}\  (\la,x,\xi) \  \text{with } \  x\notin K, \]  the map  $\si(\la,x,\xi) = p(\la,x,\xi) p(\infty,x,\xi)^{-1}$ extends to a smooth map  
\begin{equation} \label{symbo}  \sigma \colon S^q\times (R^{2n} - K\times\{0\}) \to GL(m;\C).\end{equation} 

Our bifurcation criteria will be formulated  in terms of the map $\sigma.$ In order to state our results we will need  matrix-valued differential forms. The product of two matrices of this type is defined in the usual way, with the product of coefficients given by the wedge product of forms.   The matrix of differentials  $(d\sigma_{ij})$ will be denoted by $d\sigma.$ 

 We associate to the $GL(m;\C)$-valued  function $\sigma$  of \eqref{symbo}  the  one form 
 $$ \sigma^{-1} d\sigma \  \text{defined on}\ S^q \times (\R^{2n} - K\times\{0\}).$$   Without loss of generality we can assume that $K\times\{0\}$ is contained in the unit ball $B^{2n} \subset \R^{2n}$ so that the one form  $ \sigma^{-1} d\sigma$ restricts (pullbacks)  to a   well defined one form on  $S^q\times S^{2n-1}$ which will be  denoted in the same way. Taking the trace of the 
  $(q+2n-1)$-th power of the matrix  $ \sigma^{-1} d\sigma$ we obtain an ordinary $(q+2n-1)$-form  $ tr (\sigma^{-1}d\sigma)^{q+2n-1} $ on $S^q\times S^{2n-1}.$ 

For $q$ even, we define the {\it degree}  $d(\sigma)$  of the matrix function $\si$  by 
   \begin{equation} \label{fed}
  d(\sigma)=  \displaystyle {\frac{(\frac{1}{2}q+n-1)!} {(2\pi i)^{(\frac{1}{2}q+n)}( q+2n-1)!}} \int_{S^q\times S^{2n-1}} tr (\sigma^{-1}d\sigma)^{q+2n-1}.
  \end{equation}
  
Proposition \ref{pr:fed} in Appendix C and the  integrality of  the Chern character \\  \cite[Chap. 18, Theorem 9.6]{[Hu]} 
imply that  $d(\sigma)\in \Z.$ 

 \begin{definition}\label{def:1}      
 A bifurcation point from the trivial branch  for solutions of \eqref{bvp1} is a point $ \la_* \in \R^{q}$ such that there exist a sequence  $(\la_n, u_n)\in\Lambda \times C^\infty (\bar\Omega)$ of solutions  of  \eqref{bvp1}  with $u_n \neq 0,$  $\la_n\rightarrow \la_*$ and $u_n\rightarrow 0 $ uniformly with all of its derivatives. 
\end{definition}
 
\begin{theorem}\label{th:40}
Let the problem \begin{equation}\label{bvp2}
\left\{\begin{array}{l} \cf\,(\la, x,u,\ldots ,D^{k}u) = 0,\, x \in  \Omega \\  \cg^i(\la, x,u,\ldots, D^{k_i} u)=0,\, x \in  \partial \Omega , 1 \le  i \le  r ,\end{array}\right. \end{equation}  verify   assumptions $ H_1, H_2$ and $H_3.$ 

If $q\equiv 0,4 \mod 8,$  there exists  at least one bifurcation point from the trivial branch of solutions  provided that $d(\sigma)$ is not divisible by $n(q),$
where  \begin{equation} \label{muca}  n(q)= \begin{cases} m(q/2) & \hif \, q\equiv 0 \mod 8\\  2 m(q/2) & \hif \, q\equiv 4 \mod 8 \end{cases} \end{equation} 
and $m$ is the number theoretic function defined at the end of section \ref{sec:1.2}.                                                                                  
\end{theorem} 

 Theorem \ref{th:40}  is stronger than the usual bifurcation results. Any lower order perturbation
\begin{equation}\label{bvp2'} \left\{\begin{array}{l} \cf\,(\la, x,u,\ldots ,D^{k}u) +  \cf'\,(\la, x,u,\ldots ,D^{k-1}u) = 0,\, x \in  \Omega \\  \cg'^i(\la, x,u,\ldots, D^{k_i} u) + \cg'^i(\la, x,u,\ldots, D^{k_i-1} u) =0,\, x \in  \partial \Omega , 1 \le  i \le  r ,\end{array}\right. \end{equation}
 of (1.4) with $\mathcal F'(\la,0)=0, \mathcal G'^i(\la,0)=0$  and  such that  the coefficients of the linearization of $(\cf',\cg')$ converge  uniformly  to $0$ as $\la\rightarrow \infty,$  also verifies the assumptions  $H_1$ to $H_3.$  Therefore, if $d(\sigma)$ is not divisible by $n(q),$ there must be some bifurcation point $\la\in \R^q$  for any  lower order perturbation   \ref{bvp2'} as above. 
 
\begin{remark} \label{re:5}
{\rm The definition of the degree of $\sigma$ using differential forms explains why we have assumed that $\cf$ and $\cg$  are smooth in all of its arguments including parameters.  For continuous families of linear elliptic equations with smooth coefficients, the degree of the symbol is still defined (it is called Bott's degree in \cite{[At-2]})  but it lacks of an explicit expression like the integral formula \eqref{fed}, which is due to Fedosov. One can still formulate the above theorem in terms of Bott's degree. However, its calculation in general requires a deformation of the symbol  to a simpler form.}  \end{remark}

 Now, let us consider the existence of multiple bifurcation points.
 
Putting  $\hbbr = \prod_{i=1}^r H^{k+s-k_i-1/2}(\partial \Omega; \R),$ it is shown in section \ref{sec:3.1} that,  under  the assumptions $H_1$ and  $H_2,$ the map  $(\cf , \cg)$ extends to a smooth  $q$-parameter family of Fredholm maps of index $0$ between Hardy-Sobolev spaces:
\begin{equation}\label{nbvp2}
 h=(f,g) \colon \R^q \times H^{k+s}(\Omega;\R^m ) \rightarrow   \hsb
\end{equation}
having   $\R^{q} \times\{0\}$ as a trivial branch.  Moreover, the Frechet derivative  $Dh_\la(0)$    is  the operator 
$ (L_\la, B_\la) \colon H^{k+s}(\Omega;\R^m ) \rightarrow  \hsb $ 
induced by $(\cL_{\la}, \cB_\la).$ 

Let    $\la_0 \in\Sigma(L,B)$  be an isolated singular point of $(L,B).$ We will formulate our local bifurcation result in terms of the matrix function  $R\colon S^{q-1}\rightarrow Gl(l;\R),$ where $l= \dim \Ker L_{\lambda_0},$  defined as follows: take  a small enough closed disk $D$ such that $ D \cap\Sigma(L,B)=\{\la_0\}.$ Then $R$ is defined as the restriction of  the linearization of the Lyapunov-Schmidt reduction of $h$ on a  neighborhood of $D$  to the boundary $\partial D\simeq S^{q-1}$(see  \eqref{lsbif} in section \ref{sec:3.3}).   
Since $R$ is smooth, we can consider as before the matrix  differential form $R^{-1}dR.$  
For  $ q=4s,$   we define the degree  $d(\la_{0})$ of an isolated singular point  $\la_{0} \in \Sigma (L,B)$  by: \begin{equation} \label{fedlo}
  d(\la_0)=  \displaystyle (-1)^{s+1}{\frac{(2s-1)!} {(2\pi)^{2s}(4s-1)!}} \int_{ S^{4s-1}} tr(R^{-1}dR)^{4s-1}. 
  \end{equation} 
 Much as before, by (4.19),  $d(\la_0)\in \Z.$

  \begin{theorem} \label{th:50}  Let the problem \eqref{bvp1} verify the assumptions  $ H_1, H_2$ and $H_3$ 
 of theorem \ref{th:40}.  
\begin{itemize}  
\item[i)]  If $\Sigma(L,B)$ consists only of isolated points, then they are finite in number, say
 $\{ \la_0,\dots,\la_r\},$  and  
\begin{equation} \label{multi}
   d(\sigma)  =  \sum _{i=0}^r  d(\la_i).\end{equation} 

\item[ii)]  If  $\la_0$ is an isolated singular point of $(L,B)$ and $ d(\la_0)$ is not divisible by $n(q),$  then $\la_0$ is a bifurcation point for solutions of \eqref{bvp1}.  If moreover, $ d(\la_0 ) \neq d(\sigma) \mod n(q),$ then there must be a second  bifurcation point $\la_*$ for solutions of  \eqref{bvp1} different from $ \la_0.$
   \end{itemize} 
\end{theorem}

In particular there are at least two bifurcation points  if $d(\la_0 ) \neq 0 \mod n(q) $  and either  $\sigma$ is independent from  $\la$ or $\sigma=\sigma^*$  or $ \sigma + \sigma ^* $ is  a positive definite matrix. 

  This can be seen as follows: let $S$  be any family of  pseudo-differential operators whose  principal  symbol is $\sigma.$  By \eqref{chern} and \eqref{eqfin}, $d(\sigma)$ coincides with the evaluation  of the Chern character of $\Ind S$ on the fundamental class of the sphere  $S^q.$   But in all of the above cases  the index bundle of  $S$ vanishes.  

In the first case this is clear. In the second case,  let $S'$  be  a family self-adjoint operators  with principal symbol $\sigma$ ( it is enough to take $S' = 1/2 (S  + S^*)$).   Then  $ \Ind S = \Ind S' =0,$ because  $ S'$  is homotopic  to a   family of invertible operators $S' + i \Id$  via  the homotopy $ H_t=S' + i t\Id.$  A similar homotopy  leads to the same conclusion in the third case,  using Garding's inequality. 

\begin{remark} {\rm Let us point out that, except for the one-parameter case,  the property of having isolated bifurcation points is far from being generic \cite{[Pe-4]}.}\end{remark}  
 
  \vskip 5pt
\subsection{Comments}

Our results leave many related questions open.

a) Perhaps the most interesting  one is that of {\it  global bifurcation} which predicts  the behavior of the bifurcating branch at large.  Regarding  this, the state of affairs is as follows:  the Krasnoselskij-Rabinowitz  Global Bifurcation Theorem was proved for general one-parameter families of Fredholm maps using the base-point degree  in  \cite{[Fi-Pe-Ra], [Pe-Ra]}.  Results for particular classes of Fredholm mappings of index $0$ arising from nonlinear elliptic equations and systems  are scattered around the literature. We mention  \cite{[Ki-1],[Ra-St]} among others. For a special class of  bifurcation problems involving Fredholm maps a different method was developed by  Zviagin in \cite{[Zv]}  (see also \cite{[Zv-Ra]})  using a device due to Ize. 

The extension of the Krasnosel'skii-Rabinowitz  theory  to several-parameter families of compact perturbations of identity  was  carried out mainly by the work of Alexander and Ize.  We cite  here only  \cite{[Al],[Iz-1], [Fi-Ma-Pe]} as a partial reference. The review paper \cite{[Iz-1]} has  a  wide  list of  references for this topic.  Global bifurcation for semilinear Fredholm maps was established by Bartsch in \cite{[Ba]}.  However, neither  the methods of  \cite{[Ba]}  nor the ones  in \cite{[Al]} can be used for nonlinear Fredholm maps because very little is known regarding the  extension properties of this class.  This is particularly  disappointing  since  the  bifurcation invariant  used in \cite{ [Fi-Pe-Ra]} for  the proof of the global bifurcation theorem is a particular case of our bifurcation index $\beta(f,U).$   To be precise:  taking $\Lambda =S^1$, viewed as one point compactification of the real line $\R$ and  $U=(a,b),$  under the isomorphism $\J(S^1)\equiv {\mathbb Z}_2, $  the parity  $\sigma(L,[a,b])$ used in  \cite{ [Fi-Pe-Ra]}   coincides with the local index of bifurcation points $\beta (f,U)$ considered here. 

b) Bifurcation from infinity also requires an improvement of our results.   In the case of quasilinear Fredholm maps   there is a better version of theorem  \ref{th:1.1} which, in the presence of a priori bounds, relates the order of $J(\Ind L)$ with the degree of the map $f_\la$  \cite{[Pe-1]}. This result  permits to deal at the same time with bifurcation both from $0$ and from infinity. However, the methods used here do not apply to the latter.

c) As a consequence of  the fact that our invariant depends only on the linearization of $f$ at the points of the trivial branch we have to consider not only $Bif(f)$ but  all of $\Sigma (f)$ in the formulation of the properties of the local bifurcation index. At a first glance this appears to be  an unpleasant characteristic of our  invariant since it would be preferable to deal with the set $ Bif(f) $ only.   Bartsch \cite{[Ba-1]} defined a bifurcation index of this type  for compact perturbations of  identity  parametrized by $\R^n$. It takes values in the stable homotopy group $\pi^s_n.$  In \cite{[Ba]}  his  construction was  extended to semilinear Fredholm maps. However, it is not clear how to construct  an index  of this type for general nonlinear maps. 

 On the other hand the above unpleasant characteristic  is compensated by the fact that $\beta (f,U)$ lives in $\J(\Lambda )$  which is computable  in many cases.  Indeed,  $\pi^s_n$  are still far from being completely understood  while $\J(S^n)\subset \pi^s_n$  is essentially the only known part of the stable stem.

d) As we remarked before, one of the consequences of our theory is the relation between the  nonvanishing of the Stiefel-Whitney classes of $\Ind L$  and bifurcation.   In   \cite {[Ko]}  Koschorke defined characteristic classes  of Fredholm morphisms between infinite-dimensional bundles. Koschorke's   classes are constructed  as Poincare  duals of  fundamental classes of subvarieties $\Sigma_k $ whose elements are  Fredholm operators (of index 0)  with $k$-dimensional kernel.  They are all computable from the Stiefel-Whitney classes of the index bundle. However, it is quite natural to ask whether Koschorke classes can be related to bifurcation in a  direct way.

 \vskip20pt
 
\specialsection{\footnotesize INDEX BUNDLE AND THE INDEX OF BIFURCATION POINTS}\label{sec:2}
\vskip10pt
Theorem  \ref{th:1.1},  is a special case of a  slightly more general result which is a formula relating the order of $J(\Ind L) $ in $ J(\Lambda)$  with the local multiplicity of $f_{\la}$  at $0.$  Before stating it,  we must introduce  three ingredients which appear in its formulation. 

\subsection{The index bundle}\label{sec:2.1} We shortly review  the construction of the   {\it index bundle} using a  slightly  different  approach  from the one in  \cite{[At]}  which is better suited to deal with nonlinear operators.  If $ \Lambda$ is a compact topological space,  the Grothendieck group $KO(\Lambda )$  is the group completion of the abelian semigroup $\text{Vect} (\Lambda )$ of all isomorphisms classes of real  vector bundles over  $\Lambda .$ In other words, it is the quotient of the semigroup $ \text{Vect }(\Lambda )\times  \text{Vect }(\Lambda )$ by the diagonal sub-semigroup.  The elements of $KO(\Lambda )$  are called  virtual bundles.  Each virtual bundle can be written  as a difference  $[E] - [F]$ where $E, F$ are vector bundles over  $\Lambda  $ and $[E]$ denotes the equivalence class of $(E,0).$  Moreover, one can show that $ [E] - [F]= 0$ in $KO(\Lambda )$ if and  
only if the two vector bundles  become isomorphic after the addition of a  trivial vector bundle  to both sides.  Taking complex vector bundles instead of the real ones leads to  the complex Grothendieck  group denoted by  $K(\Lambda).$    In what follows the trivial bundle  with fiber  $\Lambda\times V$ will be denoted by  $\Theta(V),$  $\Theta(\R^n)$ will be simplified to  $\Theta^n.$ 
 
    Let $X,\ Y$ be real Banach spaces and let  $L\colon  \Lambda\rightarrow \Phi(X,Y),$  be a continuous family of Fredholm operators. As before $L_{\la}\in \Phi(X,Y)$ will  denote the value of $L$ at the point ${\la} \in  \Lambda$.
Since $\coker L_\la$ is finite dimensional, using  compactness of $ \Lambda,$ one can find a finite dimensional subspace  $V \hof\,  Y$ such  that
             \begin{equation} \label{1.1}
\hbox{\rm Im}\,L_\la+ V=Y  \ \hbox{\rm for all }\  \la \in \Lambda.
\end{equation}

Because of the  transversality condition \eqref{1.1}  the family of finite dimensional subspaces 
$E_{\la}=L_{\la}^{-1}(V)$ defines a vector bundle over $ \Lambda$ with  total space  \[E= \cup_{\la \in \Lambda}\, \{\la\} \times E_\la.\] 
 Indeed, the  kernels of a family of surjective Fredholm operators form a finite dimensional vector bundle \cite{[La]}.  Denoting with  $\pi $  the canonical projection of  $Y$ onto $Y/V,$  from   \eqref{1.1} it follows that operators $ \pi {L}_{\la}$  are surjective with $ \ker\pi {L}_\la=E_\la,$ which  shows that $E\in Vect(\Lambda).$

 We define the {\it index bundle}  $\Ind L$ by: 
 
 \begin{equation} \label{defind}  \Ind L= [E]-[\Theta (V)] \in KO(\Lambda ).\end{equation}
  
If $V_1$ and $V_2$   are two subspaces verifying  the transversality condition \eqref{1.1} and $E, F$ are the 
corresponding vector bundles, we can suppose without loss of 
generality that $V_1\subset V_2$ and hence that $E$ is a subbundle 
of $F.$ The restriction of the family  $L$ to $F$ induces an isomorphism of 
$F/E$ with the trivial bundle with fiber $V_2/V_1.$  Since  exact sequences of vector bundles  split, it  follows that $F$ is isomorphic to a direct sum of $E$ with a trivial bundle  and  hence $E - \Theta(V_1)$ and $F- \Theta (V_2) $ define the same class in $KO(\Lambda ).$
This shows that  $ \Ind L$ is well defined. 

The correspondence  $ L \mapsto \Ind L $   is a natural transformation from  $\pi[-; \f]$ to $KO(-)$ which enjoys the same homotopy invariance,  additivity and logarithmic properties as the numerical index. The  proofs  of the above properties are sketched  in Appendix A. Clearly  $\Ind L =0 $ if $L$ is homotopic to a family of invertible operators.

  The index bundle of a family  of Fredholm operators of index $0,$  can be  identified with  the stable equivalence class of the vector bundle  $E$ arising in \eqref{defind}.  Let us recall that two bundles are  {\it stably equivalent}  if  they become isomorphic after  addition of trivial bundles on both sides.   Stable equivalence classes form a group isomorphic to the reduced  Grothendieck group of $\Lambda$, i.e.,  the  kernel  $ \tilde{KO}(\Lambda )$ of the rank homomorphism  $rk \colon KO(-) \rightarrow \Z. $ The isomorphism sends the equivalence class of $F$ into $ [F] - [\Theta^{r}]$ where $r=rk(F),$  \cite[Theorem 3.8]{[Hu]}.  On the other hand, the index bundle of a family of Fredholm operators of index $0$ belongs to    $\tilde{KO}(\Lambda).$

 \vskip20pt

\subsection{$\J$-homomorphism} \label{sec:2.2} 

Given a  vector bundle $E,$ let  $S[E]$ be the associated unit  sphere bundle with respect to some  chosen scalar product on $E$. Two vector bundles $E,F$ are said to be {\it stably fiberwise homotopy equivalent} if, for some $n,m,$ (and any choice of metric)  the unit sphere bundle $S ( E\oplus \Theta^n)$ is fiberwise homotopy equivalent to the unit sphere bundle $S (F\oplus \Theta^{m}).$ Let   $T(\Lambda )$ be the subgroup of $\tilde{KO}(\Lambda )$  generated by elements  $ [E]-[F] $ such that $E$ and $F$ are  stably fiberwise homotopy equivalent.  Put  $\J(\Lambda )= \tilde {KO}(\Lambda )/ T(\Lambda ).$ The projection to the quotient $\J \colon\tilde {KO}(\Lambda )\rightarrow \J(\Lambda )$ is  called  {\it  the generalized $\J$-homomorphism}.

 The group $\J(\Lambda )$  was  introduced by Atiyah in \cite{[At-1]}. He proved  that  $\J(\Lambda )$  is a finite group if $\Lambda$ is a finite $CW$-complex by showing  that  $\J(S^n)$ coincides with the image  of the stable  $j$-homomorphism of G. Whitehead (see section \ref{sec:3.3} for details).

 \vskip20pt

\subsection{Parity and topological degree}\label{sec:2.3} The third ingredient needed in order to state our main theorem is an oriented  degree theory for   $C^1$-Fredholm maps of index $0.$  The one that will be used here is the {\it base point degree}   constructed in \cite{[Pe-Ra]}.  This construction parallels the classical  approach to Brouwer degree based on regular value approximation, using an appropriate notion of  orientation for Fredholm maps.

  If $y$ is a regular value of a  proper differentiable map $f \colon\Omega\rightarrow \R^n$ defined on an open subset $\Omega$ of $\R^n,$  Brouwer's degree of $f$  on $\Omega$ is the integral number $$\deg(f,\Omega, y ) =  \sum_{ x\in f^{-1} (y)}\sgn\det Df(x). $$  In infinite dimensions  sign of the Jacobian determinant does not exists and a useful  substitute is  given by the  {\it parity } of a path of Fredholm operators of index $0$ described below.

  The  singular set $\Sigma $ of all non-invertible elements of $\Phi_0(X,Y)$ is a stratified  analytic sub-variety  of $\Phi_0(X,Y).$ Namely $\Sigma =\cup_{k\geq 1} \Sigma_k,$  where each stratum  $$\Sigma _k =\{ T\in \Phi_0(X,Y) \ / \dim \ker T=k\}$$ is  an analytic  submanifold  of $\Phi_0(X,Y)$ of codimension $k^2.$  Using  transversality, one can show that any continuous path $\ga$  in $\Phi_0(X,Y)$ can be arbitrarily approximated in norm by a smooth path $\tilde \ga $  transversal to the strata $\Sigma_k$ \cite{[Fi-Pe-0]}. By dimension counting,  a transversal path  has no intersection  with $\Sigma_k $ for $k > 1$ and only a finite number of transversal  intersection points  with the one-codimensional stratum $\Sigma _1.$  
  
  By definition,  the {\it parity} of a path $\ga$ with non-singular end points is $ \sigma(\ga) = (-1)^m$  where $m$ is the number  of intersections  with $\Sigma _1$ of  a transversal path  $\tilde \ga $ close enough to $\ga.$  It is shown in   \cite{[Fi-Pe-0]} that the parity is well defined, it is multiplicative under concatenation of paths and invariant  by homotopies  which keep end points of the path invertible. If the path is closed  its parity is defined regardless of the invertibility of the end points and  is invariant under free homotopies of closed paths.

Using the  parity  the base point degree is defined as follows.
Let   $O$ be a path connected open subset of $X.$  A $C^1$-Fredholm map  $f\colon O\rightarrow X$  is said to be orientable if for any path  $\ga $ joining two regular points of $f$  the parity of the path $Df\circ\ga$  depends only on the end points. A sufficient condition is that  $\sigma (Df\circ\ga) =1 $ for all closed paths in the domain. In particular,  all Fredholm maps of index $0$ with simply connected domain  are orientable.  
  
  Let  $f\colon O \to Y$, be an orientable Fredholm map  and let $\Omega$ be any open subset of $O$ such that the restriction of $f$ to $\Omega $ is  proper.  If the set of regular points of $f$ in $O$  is nonempty,   we choose a fixed regular point $b\in O$  (called base point) and define, for any regular value $y,$ of the map $f$ restricted to $\Omega,$ 
  \begin{equation}\deg_b(f,\Omega,y ) = \sum_{ x\in f^{-1} (y)}\epsilon(x),\end{equation}
where $\epsilon (x) = \sigma (Df \circ \gamma)$  and $\gamma$ is any path joining $b$ to $x.$ By definition, maps without regular points have degree zero.

It was proved in \cite{[Pe-Ra]} that  that this assignment extends to  an integral-valued degree theory for  proper orientable $C^1$-Fredholm maps of index $0$. The  degree is invariant under homotopies only up to sign and, as a matter of fact, there cannot be a homotopy invariant degree for general Fredholm maps extending the Leray-Schauder degree since the linear group of a  Hilbert space is connected.
However  the change in sign of  the  degree along a homotopy can be computed  using  the {\it homotopy variation property.}

An { \it admissible  homotopy}  is a continuous family of $C^1$-Fredholm maps   $$h\colon
[0,1]\times\mathcal{O}\ra Y$$  which is proper on closed bounded subsets of $[0,1]\times\mathcal{O}.$  

\begin{lemma}\label{varhomotopy}
Let     $h\colon [0,1]\times\mathcal{O}\ra Y$  be an admissible homotopy and let  $\Omega$ be  an
open bounded subset of X such that $0\not\in h([0, 1]\times\partial\Omega).$ If  $b_i\in \mathcal O$  is a  base point for  $h_i; \, i= 0,1, $  then 
\begin{equation}\label{homotopy}
\text{deg}_{b_1}(h_1,\Omega,0)=\sigma(M)\text{deg}_{b_0} (h_0,\Omega,0). \end{equation}
Here $M\colon [0,1] \ra \Phi_0(X,Y)$ is the path $ L\circ \gamma,$ where  $L(t,x)=Dh_t(x)$ and  $\gamma$ is any path in $[0,1]\times \mathcal{O}$ from  $(0,b_0)$ to $(1,b_1).$ 
 \end{lemma}

\proof  Assuming that $h$ is $C^1$ this is the content of  \cite[Theorem 5.1]{[Pe-Ra]}. In \cite{[Be-Fu]} Benevieri and Furi  used  a  very simple  argument  which allows  to  extend this  theorem to admissible homotopies  in the above sense.  We will adapt their argument to the base point degree. 

First of all we show that given a point  $t\in [0,1]$  for small enough  $\delta>0$
the homotopy property  \eqref{homotopy} holds on the interval $[t_0=t-\delta, t+\delta=t_1].$ 
 
Since $\deg_b( f,\Omega, y)$ is invariant  by small perturbations of $y,$ by Sard-Smale theorem we can assume  without loss of generality that $0$ is a regular value of $h_{t}.$ If $h_t^{-1}(0)\cap\bar \Omega $ is empty, being proper maps  closed, there exists a $\delta >0$ such that  if $|s-t| \leq \delta$ then the same holds for $h_s.$  Hence, in this case  \eqref{homotopy} is tautologically verified.   If the opposite is true, being $0$ a regular value of $h_t$ ,  $h_t^{-1}(0)=\{x_1,\dots, x_m\}.$  Applying the implicit  function theorem (in the category of continuous families of $C^1$-maps) 
on a neighborhood of each  $(t,x_i)$ and using properness we can find a $\delta>0$ such that for $s\in [t_0,t_1]$  $h_s^{-1}(0)=\{\tilde x_1(s),\dots, \tilde x_m(s)\} $ where $\tilde x_i \colon [t_0,t_1]\ra \Omega$ are continuous maps with $\tilde x_i(t) =x_i. $ Taking $\delta$ small enough we will have also  that each $\tilde x_i(s)$ is a regular point of $h_s.$ If $b_0,b_1$ are  base points for $h_{t_0}$   and $h_{t_1}$ respectively, then

\begin{equation}\label{smallhom} 
\text{deg}_{b_j} (h_{t_j},\Omega,0)=\sum_{i=1}^n \sigma\big(L\circ (t_j,\gamma_i^j)\big),
\end{equation} 
where, for $j=0,1$ and $1\leq i\leq m,$   $\gamma_i^j $ is a path in $\mathcal O$  joining $b_j$ to $\tilde x_i(t_j).$

If  $\gamma$ is any any path joining $(t_0,b_0)$ to $(t_1,b_1),$ then for each $i$ there are  two ways to reach $\big(t_1, \tilde x_i(t_1)\big)$ from $(t_0,b_0).$ One (say $\mu$) is by following  first the path $\gamma $ and then the path $ (t_1,\gamma_i^1),$ while  the second (say $\mu'$)  is to follow first the path  $(t_0,\gamma_i^0)$  and after  the path $\big(s, \tilde x_i(s)\big); t_0\leq s\leq t_1.$   Since $[t_0,t_1] \times \mathcal O $ is simply connected the two paths are homotopic and by homotopy invariance of the parity  $\sigma(L\circ \mu) = \sigma(L\circ \mu').$
But the path $ s\ra L(s, \tilde x_i(s))$ has parity one, being a path of isomorphisms. Now, the multiplicative property of the parity gives 
$$\sigma\big(L\circ  (t_0,\gamma_i^0)\big)= \sigma(L\circ \gamma)  \sigma\big(L\circ (t_1,\gamma_i^1)\big)$$  from which, taking in account \eqref{smallhom},   follows the homotopy property \eqref{homotopy} on $[t_0,t_1].$

 The general case follows again  from the multiplicative property of the parity  by subdividing $[0,1]$ in small enough subintervals.\qed

The remaining properties of a degree theory including additivity and excision hold true without change.

 \vskip20pt
 
\subsection{The main formula}\label{sec:2.4} Using the base point degree we can define the multiplicity of an isolated but not necessarily regular zero of a $C^1$-Fredholm map $ f\colon O\subset X \rightarrow Y.$  If    $x_0$ is an isolated solution of $f(x_0) =0$   its  { \it multiplicity}  is defined by   $ mult(f, x_0) = \deg_{b}(f, W, 0 ),$ where $W$ is a small enough open convex  neighborhood of $x_0$  and $b$ is any regular base point of $f$ in $W.$  Notice that the multiplicity is well defined because, being W simply connected, $f$ is orientable and all Fredholm maps are locally proper.   Moreover, the  absolute value  $|mult(f, x_0)|$ is independent from the choice of the base point.  

 Our main formula relates  the order of $\J(\Ind L)$ in  $\J(\Lambda)$   with the multiplicity  of an isolated zero at a given parameter  value. 
  
\begin {theorem} \label{th:1.3} 
Let $\Lambda$ be a finite connected $CW$-complex, $f\colon \Lambda\times O\rightarrow Y$ be  a  $C^1$-family of Fredholm maps of index $0$ and let  $L$ be the linearization of $f$ along the trivial branch.
Assume that, for  small enough $\delta,$ the only solutions of the equation  $f(\la ,x)=0$  with $ ||x|| \le \delta $  are those of the form $(\la ,0).$ If, for some (and hence all) $\nu \in \Lambda,$ the multiplicity  $k=|mult(f_\nu, 0)|\neq 0,$ then
\begin{itemize}
\item [ i)] the first Stiefel-Whitney class $w_1(\Ind L)= 0 $ 
\item [ii)]  for some $i\in \N, \ k^i \J(\Ind L)=0$  in $\J(\Lambda ).$ \end{itemize}
\end{theorem}

In particular  we have: 
\begin{corollary}\label{cor:1.4} Assume that for some $\nu  \in \Lambda, \ k =|mult(f_\nu,0)|$  is defined. If $ k  \neq 0,$  then  bifurcation arises
whenever  either  the index bundle  $\Ind L$ is  non orientable or  
 $ \J(\Ind L)\neq 0,$ and  $k=|mult(f_\nu, 0)|$ is prime to the order of $\J(\Lambda )$.
 \end{corollary}
 
Indeed,  the first assertion is clear. In order to prove the second  it
 is enough to observe that  the order of $\J(\Ind L) $ in $\J(\Lambda )$ divides the order of this finite group. Hence,  if there is no bifurcation, by the above theorem, $k$ cannot be prime to the order of $\J(\Lambda).$

If $L_\nu$ is invertible, then  the multiplicity $mult(f_\nu, 0)=\pm 1.$  Therefore, theorem \ref{th:1.1} is a special case of  the above corollary with $k=1.$

\begin{remark} {\rm A more precise invariant would be  the order  $ \J(\Ind L) $ in $\J(\Lambda).$ However,  we stated the conclusion of corollary \ref{cor:1.4}  in terms of the order of the group  $\J(\Lambda )$  since    in many important  cases (e.g., spheres)  the order of $\J(\Lambda )$ is  known.  For general parameter space without  $2$-torsion in homology  it can be estimated in terms of the homology of $\Lambda$  with coefficients in $\J(S^q)$.  On the contrary the order of  $ \J(\Ind L) $ is a rather elusive object. There is a parallel theory in terms of codegree of the index bundle (see  \cite{[Ba]}) which gives essentially the same information as the order  of  $ \J(\Ind L),$  since both numbers have the same primes on its decomposition. However, co-degree is also difficult to compute.}  
 \end{remark}

 \vskip20pt
\subsection {Proof of  the main formula} \label{sec:2.5} First  we prove   $i).$ Chose a point  $\nu\in \Lambda.$   Since $ \Lambda$ is connected,  the Hurewicz homomorphism $ h \colon\pi_1(\Lambda,\nu)\rightarrow H_1(\Lambda;\Z)$ is surjective. Therefore,  in order to show that $w_1( \Ind L)=0 $ in $H^1(\Lambda ; \Z_2)$ it  is enough to check that  $<w_1( \Ind L\circ \ga); [S^1]>= 0$ in $\Z_2 $ for any closed path $\ga \colon S^1\rightarrow  \Lambda$  with $\ga(0)=\nu =\ga(1).$  For this  we will use  the following proposition which relates the parity to the index bundle:

 \begin{proposition} \cite[Proposition 2.7] {[Fi-Pe-2]}\label{detbundle}
  Given a family $L \colon \Lambda \rightarrow \Phi_0(X,Y),$ for any  closed path $\ga \colon S^1\rightarrow  \Lambda,$
    \begin{equation} \label{detbund}  \sigma(L\circ\ga) = (-1)^\varepsilon,
  \end{equation}
 where $\varepsilon={ < w_1(\Ind  L); \ga_*([S^1])>}.$    
      \end{proposition}
 
By proposition \ref{detbundle} we have  to  show that $\sigma ( L\circ \ga)=1$  for any  closed path $\ga$ in $\Lambda$  based at $\nu.$  Let us choose a regular base point $ b\in B(0,\delta)$ for $f _\nu$ (there must be at least one since $mult(f_\nu,0) \neq 0$).  Let $L^b(t)=D_x f(\ga(t),b)= Df_{\gamma(t)}(b).$ 

Since the parity of a closed path is invariant under free homotopies,  the homotopy of closed paths  $\eta(t,s) = D_x f(\ga(t), sb), \  0\leq s\leq 1, $  shows that  
\be\label{parity0}   \si (L^b) = \si(L\circ \ga).\ee 
 
 Let  $h\colon [0,1]  \times B(0,\delta) \ra Y$ be the  homotopy defined by  $h(t,x) =f(\ga(t), x).$
  By assumption, there are no zeroes of $h$ on $ I\times \partial B(0,\delta).$ Hence, we can apply the homotopy property  \eqref{homotopy} of the base point degree to $h.$ Since $D_x h(-,b)=L^b,$ we get 
  \[ deg_b (f_\nu ,B(0,\delta))=\si(L^b)\,deg_b (f_\nu ,B(0,\delta)).\] 
 From which, being $ deg_b (f_\nu ,B(0,\delta))  \neq 0,$  we conclude  that  $\si(L\circ \ga)=\si (L^b) =1.$ 
 
  This proves the first claim.

 For the second, we will incorporate parameters into  a global version of the Lyapunov-Schmidt reduction (see section \ref{sec:3.3}) found by  Renato Caccioppoli  in \cite{[Ca]}  whose  rigorous formulation in modern terms is due to  Sapronov  \cite{[Bo-Sa-Zv]}.
 
  Let us choose an $n$-dimensional subspace $V$ of $Y$ such that the transversality condition  $\im L+ V=Y$  holds  for any  ${\la} \in\Lambda.$   Using compactness of $ \Lambda$  we can find a small enough ball $B=B(0,\delta)$ such that   the equation $f(\la ,x)$ has only trivial solutions on  $O= \Lambda\times B(0,\delta)$  and moreover
\begin{equation}
\label{eq:1.2}
\im Df_{\la}(x) + V=Y  \ \hbox{\rm for any }\  (\la ,x) \in  O.
\end{equation}
Let $\pi_{V}$ be a projector  onto the subspace $V$ and let $Z=\im (\Id -\pi_{V}).$  We split $Y$ into a direct sum   $ Y=V\oplus Z$ and we write the map  $f$ in the form $f=(g,h),$ where 
 $g \colon O \rightarrow V$  and $h \colon  O\rightarrow Z $  are  defined by $g=\pi_{V}f$ and  $h= (\Id-\pi_{V})f$ respectively.  
  
  Clearly  \eqref{eq:1.2} implies that  for each ${\la} \in\Lambda$ and  $ x \in B(0,\delta)$  the differential  $D h_\la (x) $ is surjective. Thus $h_{\la}\colon B(0,\delta) \rightarrow Z $ is a submersion for all ${\la} \in  \Lambda$ and \\ $M_\la  = h_{\la}^{-1}(0) = f_{\la}^{-1}(V$)  is a finite dimensional submanifold of $B=B(0,\delta).$ 
 
  By dimension counting, $\dim M_\la = n.$ The tangent space  to $M_{\la}$ at $0\in M_{\la}$  is  $E_\la= \ker Dh_{\la}(0) = L_{\la}^{-1}(V).$ In particular, $\Ind L $ is the  stable equivalence class of the vector bundle $E=\cup_{\la \in  \Lambda}  \,  \{\la \}\times E_{\la}.$

  Since $E$ is a finite dimensional subbundle of $\Theta(X),$ there is a family  $ \pi \colon { \Lambda} \rightarrow \mathcal{L} (X)$  of projectors with $\im \pi_\la = E_{\la}.$   We will consider $\pi $ as a vector bundle morphism from  $\Theta(X) $ onto $E.$  Let $ \phi \colon \Lambda\times B \rightarrow \Theta (Z)\oplus E $   be the (nonlinear) fiber bundle map over $ \Lambda$  defined by  $\phi(\la ,x) =(\la, h(\la ,x), \pi_{\la}(x)).$   Since $\ker Dh_{\la}(0) = \im \pi_{\la},$    $D\phi_{\la}(0)$  is  an isomorphism for each ${\la} \in\Lambda.$  

\begin{lemma}\label{le:1}  The restriction of the map  $\phi$ to a neighborhood of the zero section $  T=\Lambda\times\{0\}$  in $\Theta(X)$  is a  fiberwise differentiable homeomorphism of this neighborhood with a neighborhood of $T$ in  $\Theta (Z)\oplus E.$
\end{lemma} 

\proof  We first show that  the restriction of $\phi$ to a neighborhood of the zero section $ T=\Lambda\times\{0\}$  in $\Theta(X)$   is a fiberwise differentiable local homeomorphism using the contraction mapping principle proof of the inverse mapping theorem in  the category of spaces over a base \cite{[Cr]}.  

 Given a point $(\nu,0)\in T,$ we take a trivialization $\tau \colon E\mid_{N}  \rightarrow N\times V$of $E$  on a neighborhood $N$ of $\nu.$ Let  $ \rho \colon \Theta(Z) \oplus E \ra \Theta (Z)\oplus \Theta(V)\simeq \Theta(Y)$    be the  bundle isomorphism over $N$ defined by $\rho_{\la}(z,e)= (z,\tau_{\la}(e)).$  Composing the map $\rho \phi$ on the right with $D\phi_\nu^{-1} (0)\rho^{-1}_{\nu}$, we obtain a map  $ \bar \phi \colon N \times B \ra N\times X$ such that $ D\bar \phi_\nu(0)= \Id.$  
  
  Since $D\phi_\nu^{-1} (0)\rho^{-1}_{\nu} \rho$ is an isomorphism, we have only to prove that $\bar \phi$ is a local homeomorphism at $(\nu,0).$ In order to show this, eventually by taking smaller neighborhood of $(\nu,0)$  we can  assume that   $|| x- D\bar\phi _\la(x)||\leq \frac{1}{2} ||x||$ for all $x\in B(0,2\delta)$ and ${\la} \in N.$  Then, for each $\la$ in $N,$  the map  $\bar \phi_\la \colon B(0,\delta)\rightarrow  B(0,\frac{1}{2}\delta)$ is a homeomorphism (in fact  a $C^1$  diffeomorphism)  because  $ c_y (x )=  x-\bar \phi_{\la}(x) +y $ is a contraction on $\bar{B}(0,\delta)$ for any  $y\in  B(0,\frac{1}{2}\delta).$  We claim that the map \[\bar \phi^{-1}\colon N \times  B(0,\frac{1}{2}\delta)\rightarrow N\times B(0,\delta)\] is continuous. 

Since $D\bar \phi^{-1}_{\la}(y)$  is continuous in both variables $\la$ and $y$,  taking  $N$  and $\delta$  small enough, we have 
   $||D\bar \phi_{\la}^{-1}(y)|| \leq K $  on $N\times B(0,\frac{1}{2}\delta)$ and therefore $$ ||\bar \phi^{-1} (\la ,y) - \bar\phi^{-1} (\la ,z)|| \leq K ||y-z|| $$ there.  On the other hand,  by the continuous dependence on parameters of the fixed point  of a contraction,  $\bar \phi^{-1}(-,y)$ is continuous in the variable $\la$  for each fixed $y.$  The continuity of $\bar \phi^{-1}$ follows  from this two facts. Thus $\bar \phi $  is a local homeomorphism and hence so is $\phi$.  
   
 Finally, we observe that the restriction of  $\phi$ to the zero section $T$ is injective. It is easy to show  that if a local homeomorphism  is injective on a compact subset, then it is a homeomorphism on a neighborhood of this set. This proves the lemma. \qed\vskip5pt

Let  $U$ and $W$  be open neighborhoods   of $ T$ in $\Theta(X)$ and   $\Theta (Z)\oplus E $  respectively such that  $\phi \colon U\rightarrow W $ is a fiberwise differentiable fiber preserving  homeomorphism between them. Then the map $\psi \colon E\cap W \rightarrow \Theta (X) $ defined by  $ \psi(v) =\phi^{-1}(0,v)$ is a  fiberwise  differentiable homeomorphism  of a neighborhood of the zero section in $E$ with its image and moreover  $\psi_{\la}( W\cap E_\la ) \subset M_{\la}$ for each ${\la}  \in\Lambda.$
 
Now, we will use the map $g \colon  O \rightarrow V.$ Let $ D(E) =D(E,r) \subset E $ be a closed disk bundle of radius $r$ contained in $E\cap W$ and let $S(E)=\partial D(E)$ be the associated  sphere bundle.   Since   $\psi_\la $  sends $ D_{\la}(E) -\{0\}$ into $ M_{\la} -\{0\}$  and since $f_{\la}(x)=0$ only if $x=0,$  we have that  $|| g \psi (v)|| \neq 0$  for  any  $v\in S(E).$  Hence,  if $S^{n-1}$ is the unit  sphere in $V,$ we get a  fiber bundle map  $\bar g \colon S(E)\rightarrow  \Lambda \times S^{n-1} $  defined by
\begin{equation}\label{eq:1.3}\bar{g}(v)= \bigl( \la ,\Vert g\psi(v)\Vert^{-1} g\psi (v)\bigr).
 \end{equation}
First, we will show that the degree of the map $ \bar g _\nu  \colon S(E_\nu) \rightarrow S^{n-1}$   is  $\pm k.$
  In what follows,  if $M,N$ are oriented finite dimensional manifolds  of the same dimension,  $\Omega \subset M$  an open subset and  $f\colon {\Omega}\rightarrow N$   is  a map such that $f^{-1}(p)$ is compact,  we will denote by $\deg(f, \Omega, p)$ the Brouwer degree of $f$ in $\Omega $ with respect to $p.$  We will  use $\deg(f)$ to denote the total degree of a map between compact manifolds. 
 
 The homomorphism 
$\bar g_{\nu}^*\colon H^{n-1}(S^{n-1}; \Z)\rightarrow H^{n-1} (S(E)_\nu ,\Z)$  induced by $\bar g_{\nu}$ in singular cohomology coincides with the multiplication by $\pm \deg ((g_\nu \psi_\nu , D(E)_\nu ,0)$ (see for example \cite[Proposition 2.6 ]{[Pe-1]}). Thus $\deg(\bar g_\nu)= \pm \deg ((g_\nu \psi_\nu , D(E)_\nu ,0).$  But $\deg (\psi_\nu, D_\nu(E),0)=\pm 1 $ since  $\psi_\nu$ is a diffeomorphism.  Therefore, denoting with  $ g'_\nu$  the restriction of $g_\nu$ to  $M_\nu,$ we have 
     \[ \deg(\bar g_\nu)=\pm \deg(g_\nu \psi_\nu, D(E)_\nu,0)= \pm \deg (g'_\nu, M_\nu ,0).\] 
With the above proved, the assertion  $\deg(\bar g_\nu) =\pm k$ is a consequence  of the following reduction property of the base point degree:
 
\begin{proposition}
\label{pr:1.1} 
Let $f \colon \Omega \subset X\rightarrow Y$ be a  proper oriented $C^1$-Fredholm map of index 0. Let $V$ be an  $n-$dimensional  subspace of $Y$ transversal to $f.$ Then
$M=f^{-1}(V)$ is an n-dimensional  oriented submanifold of $\Omega.$ 
The map $g \colon M\rightarrow V$ given by the restriction of $f$ to $M$ is proper. Moreover, for any base point $b$ 
\begin{equation}\label{red}
\deg_b(f, \Omega,0)=\pm \deg (g, M, 0).
\end{equation}
\end{proposition}
\proof For $C^2$-maps  \eqref{red}  is a special case of  \cite[Theorem 5.8 ]{[Fi-Pe-Ra]}. But this theorem holds for the degree of $C^1$-Fredholm maps constructed in \cite{[Pe-Ra]}  with  exactly the same proof. \qed\vskip5pt

 From  the above proposition,  since the only zero of the map $f_\nu$ on $B(0,\delta)$ is $0,$  we get 
 $k= \deg_0(f_\nu, B(0,\delta),0) = \deg ({g'}_\nu, M_\nu ,0) = \pm \deg(\bar g_\nu),$ which proves the assertion.

 By $i),$  $E$ is an orientable  subbundle of the trivial bundle  $\Theta(X).$  Hence, we can finish the  proof of  theorem \ref{th:1.3}  using the mod-$k$ Dold's  theorem of Adams \cite[Theorem 1.1]{[Ad]}. This theorem   states that  if $E$ is an orientable vector bundle of rank $n$ over a connected  finite  $CW$-complex $ \Lambda$   and if  $\bar g\colon  S(E)\rightarrow \Theta(S^{n-1})$  is a fiber bundle map from the sphere bundle  $S(E)$ to the trivial  sphere bundle  of rank $n$ such that  for some (and hence any) $ \la \in  \Lambda$ the map  $\bar g_\la \colon  S(E_\la )\rightarrow S^{n-1}$ is of degree $\pm k$, then there exists a positive integer  $i$ such that $S(k^i E)$ is  fiberwise  homotopy  equivalent to $S(k ^i\Theta^n)$.  Thus,    $k^i \cdot  \J(\Ind L)=0$ in $\J(\Lambda )$ which proves the theorem.  \qed\vskip5pt

 \specialsection{\footnotesize A LOCAL INDEX OF BIFURCATION}\label{sec:3}
  \label{sec:locind}	
	
\vskip10pt

 Here we  will construct a local index of bifurcation points  for parametrized families of Fredholm maps. We  will consider only  families  whose range is a Kuiper space i.e.,  a Banach space with contractible linear group $GL(Y).$  Kuiper  proved that the general linear group of a Hilbert space is contractible. Later many  functional spaces   were shown to be Kuiper.  When nonempty,  the space  $GL(X,Y)$ of all isomorphisms from $X$ to $Y$ is homeomorphic to $GL(Y).$   Hence,  if  $Y$  is Kuiper, then  $GL(X,Y)$ is  contractible for any Banach space  $X$. 	
 \vskip20pt
 		 
\subsection{The local index bundle} \label{sec:3.1}

Let $L\colon U\rightarrow \Phi_0(X,Y)$ be a
continuous family of linear Fredholm operators defined on an  open set $U\subset  \Lambda$ such that $\Sigma (L)$ is a
compact subset of $U$. We define the local index bundle  as follows: 

 Let $V$ be an open  neighborhood of $\Sigma (L)$ with compact closure contained in $U.$  
Since $GL(X,Y)$ is a contractible  absolute neighborhood retract, any map from a closed subset of a metric space into  $GL(X,Y)$ can be extended to all of the space. In particular, the restriction of $L$ to the boundary $\partial V$   of $V$ can be extended to a family  $L^{\prime}\colon\Lambda - V \rightarrow GL(X,Y).$ Define
$\bar L\colon \Lambda \rightarrow  \Phi_0(X,Y)$ by patching  $L$ on $\bar V$ with
 $L^{\prime}$ on  $\Lambda-V.$ Then $\bar L$ is a family of linear Fredholm operators parametrized by $ \Lambda$
which coincides with $L$ in a neighborhood of $\Sigma (L)=\Sigma (\bar L).$

The {\it local index bundle} of the family $L$ on $U$  is defined by: 
\begin{equation}
\label{eq:2.4}
\Ind(L,U)=\Ind(\bar L)\in \tilde{KO}(\Lambda).
\end{equation}
If $V_1,V_2$ are  two neighborhoods of $\Sigma (L)$ with $V_2 \subset V_1$ (which we can always assume) 
 and if $\tilde L, {\hat L}$ are the corresponding extensions,  then $\tilde L\vert _{\partial V_1}={\hat L}\vert _{\partial V_1}$.
 
 Let  $M\colon \Lambda \rightarrow GL(X,Y)$ be
  the family defined by
\begin{equation}
\label{eq:2.5}
M_\la=\begin{cases} Id \ \qquad \hbox{\rm if}\qquad\ \la\in \bar{V_1}
\cr \tilde L_\la \circ{\hat L}_\la^{-1}\ \qquad \hbox{\rm  if} \qquad \ \la\in \Lambda - V_1. \cr
\end{cases}
\end{equation}
Then $M $  is a family of isomorphisms verifying 
$M_\la\circ { \hat  L}_\la = \tilde{L}_\la$. Since  $\Ind M=0,$  by the logarithmic property of  the index bundle,  we have 
 $\Ind(\tilde L)=\Ind(\hat L).$ This proves that the right hand 
side of \eqref{eq:2.4} is independent of the choice of $V$ and the
extension $\bar L.$ 

We will need the following additivity property of the local index bundle.

	\begin{proposition}\
	\
\label{prop:2.2}
 Let $L:U\rightarrow \Phi_0(X,Y)$ be a family  such that  $\Sigma =\Sigma (L)$ is compact. Let $U_1, U_2$ be open with $U_1\cup
U_2\subset {\it U}$ and let $\Sigma _i=\Sigma \cap U_i.$
	
	If $\Sigma _1\cap \Sigma _2 =\emptyset,\ \ \Sigma _1\cup \Sigma _2 =\Sigma$
  \noindent and if  $L_i, \ \ i=1,2 $ are the restrictions of $L$ to $U_i,$ then
	
 \begin{equation}
\label{eq:2.7}
\Ind (L,U)=\Ind(L_1,U_1)+\Ind(L_2 ,U_2).
\end{equation}

\end{proposition}	
	
\proof  Since the  index bundle  is invariant by composition with families of
isomorphisms  and since $X$ is isomorphic to $Y$ whenever  $\Phi_0(X,Y)$ is not empty, there is no loss of generality in assuming  that $X=Y$.
	
	Let $V_i$ be open neighborhoods of $\Sigma _i$ with $\bar V_i \subset U_i, 
\,  i= 1,2 ,$ and such that $\bar V_1\cap \bar V_2=\emptyset.$ Let $\bar L_i$
be  extensions of ${L_i}\mid_{\bar V_i}$ obtained as in \eqref{eq:2.4}. Using once
more the fact that $GL(Y)$ is an absolute retract we can construct two families of
isomorphisms parametrized by $ \Lambda,$ say $M_i\colon \Lambda\rightarrow GL(Y), \  i=1,2$
such that
	
	\begin{equation}
\label{eq:2.8}
\left\{\begin{array}{lll} M _i\vert \bar V_j= \Id & \text{if} & i=j \\M _i\vert \bar V_j= L_i^{-1} & \text{if}  & i\neq j.\end{array}\right.
\end{equation}

	Put $\bar L=M_2 \bar L_2\bar L_1 M _1$ and  $V=V_1 \cup V_2$. It
follows from \eqref{eq:2.8}  that $\bar L_{\la}=L_{\la}$ if ${\la} \in \bar V$ and that $\bar L_{\la}\in
GL(X,Y)$ if ${\la}  \not \in V$. By definition of the local index bundle,
\begin{eqnarray*}
 \Ind(L,U)=\Ind(\bar L)=\Ind(M_2\bar L_2\bar L_1M_1)= \\
 \Ind(\bar L_1)+\Ind(\bar L_2)=\Ind(L_1,U_1)+\Ind(L_2, U_2).
\end{eqnarray*} \qed\vskip5pt

In what follows we  will also use $ \Ind_\Lambda(L,U)$ to denote the local index bundle 
when we want  to show the dependence  of this element on  the parameter space. 

From functoriality  of the index bundle  we obtain  the following relation  between  $ \Ind_\Lambda(L,U)$
 and  the local index with respect to the one point compactification $U^+.$  
 
\begin{equation}
\label{eq:2.15}
 \Ind_\Lambda (L,U)=q^*(\Ind_{U^+}(L,U)),
\end{equation}
\noindent where $q\colon \Lambda \rightarrow U^+$ is the map collapsing  $\Lambda -U$ to the point at infinity.

 The relation \eqref{eq:2.15} suggests a  different construction of the local index bundle which works for general Banach spaces.
   This alternative approach uses $K$-theory with compact support.  We review shortly this theory below since we will need it in the sequel. 
   
     If $Z$ is a locally compact space, by definition $KO_c(Z)$ is the reduced  Grothendieck group $ \tilde{KO}(Z^+) $ of the one-point compactification $Z^+$  of the space $Z.$ However, there is  a different description of this group  in terms of  virtual bundles with compact support  \cite{[Fe],[At-Si-1]}.
 
  A {\it virtual bundle with compact support} is an equivalence class  $[E,F,a]$  of a  triple  $ (E,F,a),$  where $E,F$ are finite dimensional real vector bundles over $Z$ and where  $a\colon E\rightarrow F$ is a vector bundle morphism which is an isomorphism on the complement of a compact subset of $Z.$  Any  compact set with the above property is called { \it support}. 
  
   A triple having  an empty support is called {\it trivial}. In the set of triples  there is an obvious notion of direct sum and isomorphism. We  define an equivalence relation by  saying that  two triples $ \eta_1$ and   $\eta_2$ are equivalent provided that there are trivial triples $ \theta_1,  \theta_2$   such that $\eta_1\oplus \theta_1$  is isomorphic to $\eta_2\oplus \theta_2.$   The set of all equivalence classes is a group,  isomorphic to   $KO_c(Z).$  The complex $K$-theory with compact support  $K_c(Z)$ admits an analogous  description.

   The above isomorphism can be constructed  as follows: given a triple $(E,F,a)$ and a  relatively compact open  neighborhood $V$ of its support, by compactness  there exists a vector bundle $G$ over $\bar V$ such that $F\oplus G \cong \theta^m.$  Taking  $E'=E\oplus G$ and $a' = a \oplus \Id$ we get a triple $(E' ,\theta^m, a')$    over $\bar V$ such that  $a'$ is an isomorphism of $E'$ restricted to $\partial V$ with the trivial bundle $\partial V\times \R^m.$   We  use $a' $ in order to perform the clutching construction (see section \ref{sec:3.3}) of $E'$ with the trivial bundle  over $ Z^+ - V$ and obtain a bundle $E''$ over $ Z^+.$  It is easy  to see  that the map  $[E,F,a] \rightarrow [E''] - [\theta^m] $ is an isomorphism.  Its inverse sends  $[E]-[\theta^m] \in \tilde{KO}(Z^+)$ to the class $[E', \theta^m, a],$  where $E'$ is the restriction of $E$ to $Z,$ and $a$ is any extension to $ Z$ of a trivialization of $E$ on an open neighborhood $U$ of $\infty$ in $Z^+ $ restricted to $Z\cap U.$

Let  $Y$ be a general Banach space.  We define $\Ind(L,U)$ of a  family with compact support $L\colon U\rightarrow \Phi _0(X,Y)$ as follows:  using a finite covering of the support
 we can  find a finite dimensional subspace $F$ of $Y$ such that $\im\,L_{\la}+F=Y$ for each ${\la} \in U$. Then the family of vector spaces  $E_{\la}=L_{\la}^{-1}(F)$ is a vector bundle with a natural trivialization at infinity $a \colon E\rightarrow U\times F,$ where  $ a_\la = L_{\la}$ restricted  to $E_{\la}$.  Thus $[E, U\times F, a] $ defines an element  of $ \tilde{KO_{c}}(U)$ and it is easy to see that this element  is independent of the choice of $F$ as above. By definition,  \begin{equation}
\label{eq:2.16}
 \Ind_\Lambda (L,U)=q^*[E,U\times F,a],  
\end{equation}
 where $q$ is as in \eqref{eq:2.15}.
The relation \eqref{eq:2.15} shows that the above definition coincides with the one in \eqref{eq:2.4}  when $Y$ is a Kuiper space.
	
	 \vskip20pt
	 
\subsection{Definition and properties of $\beta(f,U)$}  \label{sec:3.2}
  
  Let $Y$ be a Kuiper space,  let $U$ be an open subset of a finite connected CW-complex $ \Lambda$ and let $O$ be an open subset of a Banach space X.
  
   Let $f:U\times O \subset X  \rightarrow Y$ be a family of $C^1$ Fredholm maps  such that $f(\la ,0)=0.$  The map  $f$ can be written in the form $ f(\la ,x)=L_\la x+g(\la ,x)$ where, as before, $ L_{\la}=Df_{\la}(0)$  and  $g(\la ,x) = o(\Vert x\Vert).$   In particular,  $Dg_{\la}(0) =0$ for all ${\la} \in  U.$
  
Recall  from  \ref{sec:1.2} that a  pair $(f,U)$ as above   is called admissible if $\Sigma(f)$  is a  compact,  proper subset of the open set  $U.$ 	
 
 The {\it local bifurcation index }  $\beta (f,U)$ of an admissible pair is defined by	
 \begin{equation}
\label{eq:2.6}
\beta (f,U)=\J(\Ind(L,U)).
\end{equation}

The rest of this subsection will be devoted to the verification of properties 
 $(\text{B}_1)$ to $(\text{B}_5)$. Property  $(\text{B}_6) $ will be proved in the next subsection. Below  we will use $\Sigma(f)$ to denote the singular set $\Sigma(L)$ of the  linearization of $f$ along the trivial branch. 

We will recast the verification of the existence property $(\text{B}_1)$  to our theorem \ref{th:1.1} by constructing  a family $\bar f \colon \Lambda \times B(0,r)\rightarrow Y $ of $C^1$-Fredholm maps  verifying the hypothesis of this  theorem and such that:
\medskip
\begin{itemize}	
\item[i)]  $ \Sigma (\bar f)=\Sigma (f) $
\item[ii)]  $\bar f$  coincides with $f$ in a neighborhood of  
	 $\Sigma (f)\times \{ 0\}$ in $\Lambda \times X$.
\end{itemize}
The construction of $\bar f$ goes as follows:  we take   an open subset $ V$
of $U$ such that  $\Sigma (f) \subset V \subset \bar{V} \subset U $. 
Arguing  as in \eqref{eq:2.4}, we  extend $L\mid_{\bar V}$  to a continuous family 
$\bar L$  defined on $\Lambda$ such that 
$\bar L_{\la}\in GL(X,Y)$ for ${\la} \in\Lambda - V$. By definition,	
	\begin{equation}
\label{eq:2.12}
 \beta (f,U)=\J(\Ind\bar L).
\end{equation}
 Let $\phi$ be   a continuous function on $ \Lambda$ with $0 \leq \phi\leq 1$,  $\phi \equiv 1$ on
 $\bar V$ and  $\phi \equiv 0$ on  $\Lambda -U.$    For $(\la ,x) \in\Lambda\times X$  we define

 \begin{equation*}
\bar g (\la ,x) = \begin{cases}
 g(\la ,\phi(\la ) x) & \mathrm{ for}  (\la ,x) \in U\times X \\ 
     0  &\mathrm{ for}  (\la ,x) \notin U\times X . \end{cases}
\end{equation*}
 Then $\bar g $ is a continuous family of  $C^1$-maps and clearly $D\bar g_\la(0) =0 $.

Finally let us  define $\bar  f$ by   $ \bar{f}(\la ,x)= \bar L_\la x+ \bar g(\la, x).$  Then $D\bar{f}_\la(0) = \bar L_\la$ and therefore, for 
 small enough  $r$, the restriction of $\bar{f}$  to  $\Lambda \times B(0,r)$ is a continuous family of  $C^1$-Fredholm maps.
Clearly the map  $\bar f$ verifies the required conditions  since it coincides  with $f$ on $ V \times B(0,r)$ and has the same singular set.

It follows from \eqref{eq:2.12} that $  \beta (f,U)=\J (\Ind \bar L) = \beta( \bar f).$
Hence, if $\beta (f,U)$ does not vanish in $\J(\Lambda ),$ by theorem \ref{th:1.1},  there must be
a bifurcation point of $\bar{f}$ belonging to  $\Sigma(\bar L) = \Sigma(L).$  Since $f$ coincides with $\bar{f}$
 on $V\times B(0,r),$  this point must be a bifurcation point for
$f$ as well. This  completes the verification of $(\text{B}_1).$ 
\smallskip

That $(\text{B}_2)$ holds is clear from the definition of the local index bundle. In order to prove  the additivity property $ (\text{B}_4),$  it is enough to consider  the case of two open sets. Notice that, being $\Lambda$ connected,  if $(f,U)$ is admissible so are $(f_i,U_i).$ Then $(\text{B}_4)$  follows from proposition \ref{prop:2.2}  applying the functor $\J $ to both sides. 

 In order to show $(\text{B}_5)$ let us notice that, if 
 $\alpha \colon Q\rightarrow  \Lambda$ is continuous, then by 
functoriality of the index bundle	
	\begin{equation}
\label{eq:2.10}
 \Ind(L\circ \alpha, \ \alpha ^{-1}(U))=\alpha^* \Ind(L,U).
\end{equation}
		If $g=f\circ (\alpha \times id_X),$ then $Dg_\la(0)=L\circ\alpha(\la).$  Applying 
the functor $\J$ to \eqref{eq:2.10} we obtain $(\text{B}_5).$ 	

The homotopy invariance property $ (\text{B}_3)$  follows from  $(\text{B}_5),$  
since an admissible  homotopy $h$ is nothing but an admissible family of  $C^1$-Fredholm maps  parametrized by  the open subspace   $V=[0,1] \times U$  of the space $\Gamma= [0,1]\times  \Lambda.$ Thus $h$  defines  an
element $\beta (h,V)\in \J(\Gamma)$. By $ (\text{B}_5),$ denoting with  $ i_0$ and $i_1$  are the  top and bottom inclusion of $\Lambda$ in $\Gamma,$    we have $ \beta (h_j,U)=i_j^* \beta(h,V),\,  j=0,1.$ But $i_1^*= i_0^*$ and hence $\beta (h_0,U)=\beta (h_1,U).$

\br{\rm  If $Y$ is a general Banach space, and we put  $\beta(f,U)= \J(\Ind (L,U)),$ where $\Ind (L,U)$ is defined  by \eqref{eq:2.16}, then we can show that $ \beta(f,U)$ verifies properties  $(\text{B}_2)$ through $\text{B}_6$ using  standard properties of  $\tilde{KO}$ as  generalized cohomology theory.  However, the crucial property  $(\text{B}_1)$ is missed in this setting because our proof of  $(\text{B}_1)$  relies on an extension property which does not hold for general Banach spaces.}\er 

\subsection{Comparison with the Alexander-Ize  invariant} \label{sec:3.3}
 
In order to complete the proof of  theorem \ref{th:locind}   we have to verify   the property   $(\text{B}_6).$
  
   We begin by introducing the Alexander-Ize invariant. Only the stable version of this invariant constructed in \cite{[Al]} will be considered here.

 Let  $f\colon \R^k\times \R^n\rightarrow\R^n$ be  a  $C^1$-family of maps from $\R^n$ to itself, parametrized by
$\R^k$. Assume that $f_{\la}(0)=0$
and let $L_{\la}$ be the derivative of~$f_{\la}$ at~$0$.
Let ${\la} _0$ be an isolated point in the set $\Sigma(L)$. The homotopy class of the restriction 
  of the map ${\la} \mapsto L_{\la}$ to the boundary of a small closed disk $D=D(\la_0,\epsilon)$ centered at ${\la} _0$ defines  an element $\gamma^{n}_f$ in the homotopy
group ~$\pi_{k-1}GL(n;\R)$ (here and below  we are using  the fact
 that  our target spaces are H-spaces and hence the free and pointed homotopy classes are the same). 

 Stabilizing $\gamma^{n}_f$ through the natural  inclusion of $GL(n)=Gl(n;\R)$ into 
$GL(n+1)$ one gets an element  $\gamma_f$ belonging to the  homotopy group
$\pi_{k-1}GL(\infty),$ where  the space $GL(\infty)= \cup_{n\geq 1} GL(n)$  is endowed with the inductive topology.  The element $\gamma_f$  is the Alexander-Ize invariant. 

Let $\pi^s_{k-1}=\text{limdir}\ \pi_{m+k-1}S^m$ be the $(k-1)$-stable homotopy stem. In \cite{[Al]} Alexander proved  that ${\la} _0$ is a bifurcation point of $f$ provided the image of $\gamma_f$ by the stable $j$-homomorphism  $ j\colon \pi_{k-1}GL(\infty)\rightarrow \pi^s_{k-1}$ does not vanish. 

The above definition  can be easily extended  to continuous  families  of $C^1$-Fredholm maps   $f\colon \R^k\times X\rightarrow Y.$ 

Indeed, assume  that $\la_0$ is isolated in $\Sigma(L)$  and let  $D $ be as  before.  A  {\it regular parametrix}  \cite{[Fi-Pe-0]}  for the family $L$ is a family of isomorphisms $A\colon D\rightarrow GL(Y,X)$  such that $L_{\la}A_{\la}=Id_{Y}+K_{\la},$  with  $\im K_{\la}$  contained in a fixed finite dimensional subspace $F$ of $Y.$ Since  $D$ is contractible, any  family  $L$ as above possesses a regular parametrix  (see the proof of  lemma \ref{lem:bell} below). 
  
Putting $ N_{\la}=(Id+K_{\la})\mid_{F},$ the map $N$  sends $\partial D$ into $GL(F).$ Choosing  a basis of $F$ we get a family of matrices in $ GL(m)$ parametrized by $ \partial D \simeq  S^{k-1}.$  By the preceding discussion, the stable homotopy class of $N_{\mid \partial D} $  defines an element 
$\gamma_f \in \pi_{n-1}GL(\infty)$
which can be  shown to be independent from the  choice of $D$ and the parametrix $A.$ 
By definition, the element  $\gamma_f$ constructed above is the {\it Alexander-Ize  invariant} of $f$ at $\la_0.$
 
Let us discuss the identification of $\J(S^k)$ with the image of the  stable  $j$-homomor\-phism of G.Whitehead.  
 
 A spherical  fibration is a  fibration locally fibre homotopy equivalent to a product of the base with an $n$-sphere.  Recall that the reduced group $\tilde{KO}(\Lambda)$ can be identified with the group of stable equivalence classes  of vector bundles over $\Lambda.$  In a similar way one can introduce the group  $\tilde{KF}(\Lambda)$ of stable fibre homotopy classes of spherical  fibrations \cite{[Di]}.  $\tilde{KF}(\Lambda)$ becomes a group under  the operation of fiberwise smash product.
 
  As in the case of $\tilde{KO},$ the group  $\tilde{KF}$ is a homotopy  functor represented  by the classifying space $BH(\infty)$ of the monoid $H(\infty) =  \cup_{n\geq 1} H(n),$  where  $H(n)$ is the space of all homotopy equivalences of $S^n.$ 
  
 Since  working directly with $GL(n)$ instead of the homotopy equivalent group $O(n)$  simplifies many  arguments in this section,  we deviate slightly from the usual  convention.  The  later defines  $H(n)$  to be the set of homotopy equivalences from $S^{n-1}$ into itself  and identifies  the previously defined  $J$ homomorphism  with the natural transformation which assigns to each vector bundle $E$  its unit sphere bundle $S(E).$ 
 
   Here instead,  we  will consider    $J \colon\tilde{KO}\rightarrow \tilde{KF}$ to be defined by the inclusion of the total space of a vector bundle in its fiberwise one-point compactification.  Since the fiberwise one point compactification of a vector bundle is  a suspension of its  unit sphere bundle,  we obtain a factorization:

  \begin{equation}\label{eq:J} 
 \bfig
\Vtriangle/>`>`<-/<500,500>[ \tilde{KO}(\Lambda)` \tilde{KF}(\Lambda)`\J(\Lambda);\J`\J`\text{inclusion}]
 \efig
 \end{equation}
which leads to the identification of $J(\Lambda)$ with the image of the horizontal arrow.

 Taking $\Lambda =S^k$ we obtain  a
   commutative diagram
 \begin{equation}\label{eq:j}\bfig
\square<800,500>[\tilde{KO}(S^{k})`\tilde{KF}(S^{k}) `\pi_{k-1}GL(\infty)`\pi_{k-1}H(\infty);\J`\partial_0`\partial_1`j]
 \efig\end{equation}

In the above diagram $j$ is the homomorphism induced in homotopy by the map which assigns
 to each element of  $GL(n)$  the obvious  extension to a map from the one point compactification $S^{n}$ of $\R^n $ into itself. The vertical arrow $\partial_0$ takes the stable equivalence class of a vector bundle $E$ over $S^k$  to the stable homotopy class of  $$\psi_T=T_{-} T^{-1}_{+} \colon S^{k-1}\rightarrow GL(n),$$  where $T _{\pm}$ are  trivializations for the restrictions of $E$ to the upper and lower  hemisphere of $S^{k}$  respectively.   The vertical arrow $\partial_1$ is defined in a similar way. The homomorphisms   $\partial_i, \, i=0,1$ are   isomorphisms whose  inverses are given by the clutching construction.  
 
 Under the identification of $ \pi_{k-1}H(\infty)$  with  $\pi_{k-1}^{s}$ via the  isomorphism established in     lemma $1.3$ of \cite {[At-1]}, the homomorphism  $j$ in \eqref{eq:j} coincides with the stable  $j$-homomorphism of G. Whitehead and
  the vertical arrow  $ \partial_{1}$  sends $ \J (S^{n})$ isomorphically onto $\im j.$  In what follows we will  identify the group $ \J (S^{n})$ with $\im j$ by means of the restriction of $\partial_{1}$  to $ \J(S^{n}).$

Before going to the verification of  $\text{B}_6$ we will need the analog of $\partial_{0}$ at the operators level.  Let   $ \partial \colon \pi_k \Phi_0(X,Y)\rightarrow \pi_{k-1} GL(\infty),$  be defined as follows: let
$ L \colon S^k\rightarrow \Phi_0(X,Y)$ be a  family representing the homotopy class $\alpha\in  \pi_k \Phi_0(X,Y).$ We can take parametrices  $A_\pm$  of  $ L_\pm =L \mid_{D_\pm} $ such that,  for any $\la \in D_{\pm},$ the operators $ K_{\pm \la }  = L_{\pm \la} A_{\pm \la} - Id $  take values in the same  $r$-dimensional subspace $F$ of $Y.$  Then, for  $\la \in S^{k-1},$ the operator  $ A^{-1}_{-\la } A_{_{+\la}}$ sends  $F$ into itself.  By definition, $\partial(\alpha) $ is the stable homotopy class of  $\phi_{A}\colon S^{k-1}\rightarrow GL(F)\cong GL(r)$ defined by  
$\phi_A (\la) =  A_{- \la}^{-1} A_{_{+\la}}\mid_{F}.$

\begin{lemma}\label{lem:bell}
 The diagram 
\begin{equation}\label{eq:gf}  
\bfig
\Ctriangle/<-`>`>/<600,400>[\tilde{KO}(S^{k})` \pi_k \Phi_0(X,Y)`  \pi_{k-1}GL(\infty);\Ind`\partial_0`\partial]
 \efig
 \end{equation}
is commutative \end{lemma}
\proof Let $F$ be any subspace of  $Y$ verifying the transversality condition \eqref{1.1}. Then the index bundle of $L$ is the stable class of  $E= L^{-1}(F).$ 
Given  trivializations $T_{\pm}\colon E\mid_{D_{\pm}}\rightarrow D_{\pm}\times F $ we construct
the  parametrices $A_{\pm}$ of $L_\pm$ as follows:
for $\la \in D_\pm$ we put 
 \begin{equation}\label{parametrix}
A_{\pm \la}= \bigl(Q'L_{ \la} + T_{\pm \la}Q_{\la}\bigr),^{-1} 
\end{equation}
 where $Q_{\la}$ is a continuous family of projectors of $X$ with $\im Q_{\la}=E_{\la}$ and $Q'$ is  a projector  with $ \ker Q' = F.$ 
 
$Q'L_{\la} + T_{\pm \la}Q_{\la}$ are injective Fredholm operators of index $0$ and hence are invertible for any $\la\in D_{\pm}.$  Thus $A_{\pm \la}$ are well defined. Moreover,  the image of $ L_{\pm \la} A_{\pm \la} - Id_{Y}$ is contained in  $F.$  Using  parametrices  $A_{\pm }$ in the definition of the homomorphism $\partial_0$ one easily checks that on $S^{k-1},$  $A_{-\la}^{-1}A_{_{+\la}}= \Id_{Y} -  K_{\la},$
 where  $K$  is such that   $\im\, K_{\la} $ is contained in
$F$. Since  $Q'L_{ \la}$ vanishes on $E_{\la},$ it follows that on $F$ the operator
 $(Q'L_{\la} + T_{+\la})^{-1}$ coincides with $T_{+\la}^{-1}$ and hence  $A^{-1}_{-\la}A_{_{+\la}}$
restricted to $F$ is nothing but  $T_{_{-\la}}T^{-1}_{+\la}$. Thus, with the above choice of parametrix, we have  $ \psi_{T} =\phi_{A}$  and therefore 
$\partial_{0}\circ  \Ind =\partial.$ \qed\vskip5pt

 \begin{proposition}
\label{pr:alex}
Let  \ ${\la} _0$ be the only  singular point of  $f\colon U\times O \rightarrow Y.$  Assume that $U\cong \R^{k}.$ Then, on $ U^{+ }\cong  S^{k},$ the identification  $\partial_{1} \colon \J(S^k)\simeq  \im\,j$ sends  $\beta_{S^{k}}(f,U)$ into  $j(\gamma_f).$ \end{proposition}
 
\proof  
 We can  assume without loss of generality that ${\la} _0$ is the
north pole of $S^k$ and take in the definition of $\gamma_{f}$ the upper  hemisphere $D_{+}$ as the disk  $D.$
 Let $L$ be the linearization of the family $f$ along  the trivial branch and let $\bar L$ be any extension of $L\mid_{D_{+}}$  to all of $S^k$ such that  
$\bar L_{\la}$ is an isomorphism for $\la\in D_{-}.$ 

Putting together  the commutative diagrams  \eqref{eq:gf} and \eqref{eq:j} 
we obtain going up and  right  $\J(\Ind \, \bar L) $ which by definition is $ \beta(f,U).$ On the other hand,  going  down and right  we get $j(\gamma_f).$  Indeed, if  $A_{+}\colon D_{+}\rightarrow GL(Y,X)$  is  the parametrix for $\bar L_{+},$ we can take 
 $A_{ -}= \bar L_{-}^{-1}.$  Then $\displaystyle{ A_{- \la}^{-1} A_{_{+\la}}\mid_{F}=  \bar L_{-\la}A_{+\la}\mid_{F}=N_{\la}},$ where  $N$ is the family defining the class   $\gamma_{f} $. \qed\vskip5pt
  
The above proposition shows  that  $\text{B}_6$ holds true in the case $\Lambda =S^{k}$. The general case now follows from this and \eqref{eq:2.15}. This completes the proof of  theorem \ref{th:locind}.

\vskip 10pt

   When the family $L$ behaves in a regular way close  to $\la_0,$  from  the  above proposition, using  results of Alexander and York in \cite{[Al-Yo]}, we  can obtain sufficient conditions for the nonvanishing of $\beta_{S^{k}} (f,U)\neq 0$ in a small enough neighborhood $U$ of $\la_0$  in terms of the dimension  of $\ker  L_{\la_0}.$ 
   
\begin{corollary}
\label{cor:2.1}
Let 
$f:\R^{k}\times X\rightarrow Y$ be a continuous family of $C^{1}$-Fredholm maps and let 
${\la} _0$ be such that  for $L=Df_{\la}(0)$ the following
condition holds: 
there exists a positive number $r$  such that for  small enough $\delta$  
\begin{equation} \label{regul} \|L_{\la}x \| \geq r \|\la -\la_{0}\| \|x\| 
\hfor 0 \le \|{\la}-\la_{0}\|\leq \delta.\end{equation}

 Let $c_k$ be defined by  \begin{equation} \label{cn}  \begin{matrix} k=&1&2&3&4&5&6&7&8 \cr c_k= &1&2&4&4&8&8&8&8 \cr\end{matrix} \qquad c_{k+8}= c_k, \end{equation}
\noindent then, for $k\equiv 1,2,4,8$ mod-$8,$   dim $\ker L_{\la_{0}} = m $ is a multiple of  $c_k.$
Moreover, if $m=dc_k$ with $d$ an odd integer, then $\beta_{S^{k}} (f, D(\la_{0}, \delta))\neq 0$ in $\J(S^k).$
\end{corollary} 

 This follows from proposition \ref{pr:alex} and  computation of $j(\gamma_f)$ in \cite{[Al-Yo]}.

\begin{remark}\label{id}
{\rm The {\it intrinsic derivative} of a smooth family $L\colon \R^k \rightarrow \fz$  at $\la\in \R^k$ is the map \[ \dot{I}L(\la) \colon \R^k \rightarrow \cL(\ker L_\la;\coker L_\la)\] defined as follows: $ \dot{I}L(\la)v$ is the restriction to $\ker L_\la$ of the  ordinary Frechet derivative  $DL(\la)v$  followed by the projection to the $\coker L_\la.$  When the family $f$ is smooth in all variables,  condition \eqref{regul} can be checked  from the intrinsic derivative of $L$ at $\la_0.$  It  was shown in \cite{[Fi]} that for smooth families the regularity condition  \eqref{regul} holds if  and only if  for any $v\in \R^k$  $\dot{I}L(\la)v$ is an isomorphism.} 
\end{remark}

In the final part of the section we will point out  the relation of our construction of $\gamma_{f}$  with the Lyapunov-Schmidt reduction. We will use this relation  in the proof of theorem \ref{th:50}.  Moreover, we will be able  to compare the approach we have chosen here with that of Ize in  \cite{[Iz],[Iz-1]}  which uses as the unstable version of  $\gamma_f $  the homotopy class of the linearization at $0$ of the reduced map. 

For simplicity, let us assume that  the isolated  singular point is $\la_0= 0.$

 Let $Q'$ and Q  be projectors on $ Y_1=\im L_0$ and $ E_0= \ker L_0 $ respectively.  Then  $F_0 =\ker Q' \simeq \coker L_0.$  Under the splitting of  both $Y$ and $X$ into a direct sum $Y_1 \oplus F_0$ and $X_1 \oplus E_0$  the Frechet  derivative  $D_{x_1} Q'f (0,0)$ in the direction of $X_1$  is an isomorphism. 

 By the implicit function theorem,  there exist a map $\rho$ defined on a neighborhood of $(0,0)$ in $\R^k\times E_0$ with values in $X_1$ such that, close enough to $(0,0)\in R^k\times X,$ we have  $Q'f(\la, x_1+ v)=0$ if and only if $x_1=\rho(\la,v).$  It follows that, for small $(\la,x),$ the solutions of $f(\la,x)=0$ are in one to one correspondence with the solutions of the finite dimensional reduced system $r(\la,v)=0$ (called bifurcation equation),  where the map  $r$  is defined on a product  neighborhood of $(0,0)$ in $\R^k\times E_0$   by 
\begin{equation} \label{bifls} r(\la,v)=(\Id-Q')f(\la,\rho(\la,v) +v).\end{equation}

Clearly $r(\la,0)=0.$ Let $  R_\la =D\,r_\la(0)$ be the linearization of $r$ at the trivial branch. Taking a small enough closed  disk $D=D(0,\delta)$  centered at $0,$ the restriction  of $R$ to $\partial D$ defines a map 
$ R \colon S^{k-1} \rightarrow  GL(E_0,F_0)$  and hence (after a choice of basis of both spaces) a family of nonsingular matrices   \begin{equation} \label{lsbif} R \colon S^{k-1} \rightarrow Gl(m), \ m=\dim E_0, \end{equation} whose homotopy class depends only on the choice of orientations of $E_0$ and $F_0.$ Let us remark  that  the bifurcation invariant defined by Ize in \cite{[Iz]} is the image of the homotopy class of $R$ by the unstable $\J$-homomorphism.  

\begin{proposition}\label{prop:ls}  With an appropriate choice of orientations the stable homotopy class of $R$  in $ \pi_{k-1} GL(\infty)$ coincides with the Alexander-Ize invariant $\gamma_f.$ 
\end{proposition}

\proof  We will show that  $R$ is homotopic to the family of matrices  $N$ used in the definition of $\gamma_f.$ This
 will prove the proposition.

Let  $S\colon Y_1\rightarrow X_1$ be the inverse of the operator $L_0$  restricted to $X_1$. An easy calculation (see \cite{[Iz]}) gives  $R_\la = (\Id -Q')L_\la M_\la,$ where $M_\la \in GL(X)$  is  defined  by   $M_\la = [\Id +SQ'(L_\la -L_0)]^{-1}.$ For small enough $D$ the transversality condition \eqref{1.1} is verified  with $F= F_0.$ Thus the family of subspaces $E_\la = L^{-1}_\la(F_0)$ form  a trivial vector bundle over $D.$ 
 
 Given a trivialization   $T \colon E \rightarrow D\times F_0, $  denoting with  $Q_\la $ the family of projectors on $E_\la,$    the family of isomorphisms
$A _\la= \bigl(Q'L_{ \la} + T_{ \la}Q_{\la}\bigr)^{-1}$ 
is a parametrix  $A$ of $L_{\mid D}.$ 
Thus, each $A_\la $ is an isomorphism and  we have $L_\la A_\la = \Id  + K_\la $  with   $\im K_\la \subset F_0$  for all $\la\in D.$ Arguing as in the proof of  lemma \ref{lem:bell} we obtain  ${A_\la}_{\mid{F_0}}= T^{-1}_\la \colon F_0\rightarrow E_\la.$  Using this in the definition of $K_\la$  we get 
  \begin{equation}\label{paramet}
 N_\la = (\Id + K_\la)_{\mid{F_0}}  = (\Id -Q')L_\la T^{-1}_\la.
\end{equation}
We write $N_\la $  in the form 
  \begin{equation}\label{parameto}
N_\la =(\Id -Q')L_\la M_\la (M^{-1}_\la T^{-1}_\la). 
\end{equation}
Observing that  $M^{-1}_\la=\Id +SQ'(L_\la -L_0)$  sends isomorphically  $E_\la$  into $E_0$ 
we have that $ H_\la = M^{-1}_\la T^{-1}_\la $ sends $F_0$ isomorphically into $E_0$  for all $\la \in D.$ Restricting our families to $\partial D$ we obtain 
$N_\la =R_\la H_\la $ and hence $N$ is homotopic to $RH_0$ via the homotopy $h(t,\la)=R_\la H_{t \la}$.  Choosing basis in $E_0$ and $F_0$  such that the determinant  of the  matrix of $ H_0$  is $1$ we obtain a homotopy  between  the matrix families $R_{\mid \partial D}$ and $N_{\mid \partial D}$.   This proves the proposition.

 \vskip20pt

\specialsection{ \footnotesize BIFURCATION OF SOLUTIONS OF NONLINEAR ELLIPTIC  BVP}\label{sec:4}

 Using results from the previous chapters  we will prove the criteria for bifurcation of nontrivial solutions of elliptic boundary value problems stated in section $1.$
 
 Our strategy will be as follows: extending  the  Agranovich reduction \cite{[Ag]} to parametrized families of elliptic boundary value problems we will show that $\Ind L$  coincides with  the index bundle of a parametrized family $\cS$  of pseudo-differential operators of order zero on $\R^n$ belonging to a class of  introduced by Seeley in \cite{[Se]}. Then we will  use  the Atiyah-Singer theorem  for operators in this class  which states that  $\Ind L$ (i.e. analytical index of the family) can be obtained  from  the symbol class by a homomorphism called {\it topological index}.  In our special case the  topological index  is an isomorphism which  coincides up to sign with the inverse of the Bott isomorphism. This  makes all  calculations simpler.  Using Fedosov's formula  for Chern character of the index bundle  and applying well  known results about the kernel of   $\J$-homomorphism, due to Adams and others,  we will obtain  criteria for nonvanishing of $\J(\Ind L)$  and hence for the appearance of nontrivial solutions of the problem.

 \vskip20pt
\subsection{The Agranovich reduction}  \label{sec:4.1}
 We will consider particular families of boundary value problems for which the reduction in the title can be carried out. We will work out the reduction for families continuously parametrized by general compact spaces  since we will need this generality in \cite{[Pe]}.  
Let  
 \[\begin{cases}
{\mathcal L}_\la(x,D)=&\sum_{|\alpha|\leq k}a_{\alpha }(\la, x) D^{\alpha },\\ 
\mathcal{B}^i_\la(x,D)=& \sum_{|\alpha|\leq k_i}b^i_\alpha (\la,x)D^{\alpha},\, 1\leq i\leq r, \end{cases}\] 
be a family of linear boundary value problems where the matrix functions \\ $a_{\alpha }(\la,x) \in \C^{m\times m}, \, b^i_\alpha(\la,x)  \in \C^{1\times m}$  are smooth in $x$  and continuously depending on a parameter  $\la $  belonging  to a compact topological space $\Lambda.$

 The class  under consideration is described  by  axioms $A_1$ to $A_3$ below. 

\begin{itemize}

\item[$A_1)$]  
 For all $\la \in \Lambda,$ the boundary value problem  $({\mathcal L}_\la(x,D),\mathcal{B}_\la(x,D))$ is elliptic.  Namely,   ${\mathcal L}_\la(x,D)$ is  elliptic, properly elliptic at the boundary, and the rows of the  boundary operator $\mathcal B_\la(x,D)=(\mathcal B^1_\la(x,D),...,  \mathcal B^r_\la(x,D))^t$  verify the Shapiro-Lopatinskij  condition with respect to $ {\mathcal L}_\la(x,D)$ (Appendix B).

\item [$A_2)$]   There exists a   $\nu \in \Lambda $ such that for every  $f \in C^\infty(\bar \Omega; \C^m)$ and $g \in C^\infty (\partial\Omega; \C^r)$  the problem 
\[\left\{ \begin{array}{l}{\mathcal L}_\nu (x,D) u(x) = f(x)\ \text{for}\  x \in \Omega \\ 
\mathcal {B}_\nu (x,D) u(x) = g(x)\ \text{for}\  x\in \partial\Omega
\end{array}\right. \]
has  a unique smooth solution.
\item [$A_3)$]  \begin{itemize}
\item[ i)]  The coefficients  $b^{i}_\alpha(\la,x),  \  |\alpha|=k_i,  1\leq i\leq r ,$
 of leading terms of boundary operators $\mathcal B^1_\la(x,D),\dots, \mathcal B^r_\la(x,D)$ 
 are independent of  $\la.$  

 \item[ ii)]  There exist a compact set  $K\subset \Omega $ such that  the coefficients \\$a_{\alpha }(\la, x), \, |\alpha |=k $ of  leading terms of $\cLL,$ are independent of $\la $ for $x\in \bar\Omega-K.$ \end{itemize}\end{itemize}

Under assumption  $A_1$   the differential  operators   $(\cL_{\la}, \cB_{\la})$  define  a continuous family of  bounded  semi-Fredholm operators  (Appendix B) 
\begin{equation}\label{fred}
( L,B) \colon  \Lambda  \rightarrow \mathcal{ L} (H^{2k+s}(\Omega; \C^m );  H^{s}(\Omega; \C^m )\times\hbb).  
\end{equation}  

  By   $A_2$ and  the regularity of solutions of elliptic equations,  the kernel of the  operator   $(L_\nu, B_\nu)$ reduces to $u\equiv 0$  and its  image contains a dense subspace. Therefore, $(L_{\nu}, B_\nu)$   is an isomorphism which on its turn, by  the invariance property of the index,  shows that  the family  $(L,B)$  is a continuous family of Fredholm operators of index $0.$ 
  
   We will show  that the index bundle of the family $(L,B)$ coincides  with  the index bundle of a  family of  a particular class of pseudo-differential operators on $\R^{n}$ introduced by Seeley in \cite{[Se]}. 
  
A {\it symbol of class $S^{k}(\co)$}  is a function 
$ \rho \in C^ \infty (\co\times \R^n ;\C^{m\times m})$ defined on an open subset $\co$ of  $\R^n,$ verifying  following property:

 for  every compact   subset  $K$ of $\co$ there is a constant $C$ such that,  for $x \in K,$ 
\begin{equation}\label{boundpd}
| D^\alpha _x D^\beta_\xi \rho (x,\xi) | \leq C( 1+ |\xi|^{k-\beta}).
\end{equation} 
The set $S^{k}(\co)$  is naturally  a Frechet space with the topology induced  by the family of seminorms 
\begin{equation}\label{frect}
\pi_{k,K} ^{\alpha \beta}(\rho)= \displaystyle{ sup_{ x\in K, \xi\in \R^n} } (1+ |\xi|) ^{\beta-k}|D^\alpha _x D^\beta_\xi \rho (x,\xi)|.
\end{equation}
A  {\it  pseudo-differential operator of order $k$ } acting on the space  $\cd (\co)^m$ of all  smooth  $\C^m$-valued  functions  $u$ with compact support in $\co$   is  defined  by an integral
 
\begin{equation}\label{pdo}
\cq u(x) = (2\pi)^{-n} \int _{\R^n} e^{ix\xi} \rho(x,\xi) \hat{u}(\xi) \, d\xi,
\end{equation}
where  $\rho \in S^{k}(\co)$ and  $\hat{u}$ denotes the Fourier transform of  $u.$     
Every  pseudo-differential operator $\cq$ of order $k$ extends  to  a
linear continuous map 
\[Q\colon \hck\rightarrow \hl.\]
 Here   $\hl,$   is the space of  $\C^m$-valued distributions  $u$ on  $\co$  such that,  for all $\varphi \in \cd(\co),  \, \varphi u \in \hsn,$   with the topology induced by the family of semi-norms $||\varphi u||_s.$  The space  $\hck $ is the union over  all compact subsets $K$ of $\co$ of  $$H^{k+s}_K (\co; \C^m) =\{ u\in  H^{k+s}_{loc}(\co; \C^m)|\, supp\, u \subset K\} $$  endowed with the
  direct limit topology  for the family of inclusions.  
\vskip 5pt
A pseudo-differential operator $\cL$  of order $k$  is  said to be elliptic  if it possesses a  (rough) {\it parametrix}  or {\it regularizator.} This  is a proper(\cite{[Ch-Pi]})  pseudo-differential operator  $\cP$   of order  $-k$ such that  both  $\cL \circ\cP -\Id $ and  $\cP\circ\cL -\Id$ are of order $-1.$    A  stronger  notion  of parametrix is used in regularity theory  but  for the  purpose of computing the index bundle  this one  will be sufficient.  

Elliptic differential operators are elliptic in the above sense.  As a parametrix  of  $\cL$   one can  take  the pseudo-differential operator $\cP$  associated to the  symbol

\begin{equation}
\label{par}
 \rho (x,\xi) = \phi (|\xi|)p^{-1}(x,\xi),\, \text{if } \, x\in \co,\\
\end{equation}  

where  $p= \displaystyle \sum_{|\alpha|= k}  a_\alpha(x)\xi^\alpha$ is the principal symbol of $\cL$ and $\phi$ is a smooth function with $ \phi (r) \equiv 1 $ for $ r\geq 1$ and  $\phi (r) \equiv 0$ on a  small  neighborhood of $0$.

We will deal only  with pseudo-differential operators whose symbols  enjoy a further property :

 Outside of a small neighborhood of $\co \times \{0\}$ 
\begin{equation}
\label{symbol}
\rho(x,\xi)= \rho_s(x,\xi) + \delta(x,\xi),\,\text{where}\, \rho_s = \lim_{\mu\rightarrow \infty} \rho (x , \mu \xi)\mu^{-s} 
\end{equation}
is a homogeneous function of degree $s$ defined on $ \co\times(\R^n-\{0\})$ and $\delta $ is a symbol of order $s-1.$

  This class of  pseudo-differential operators contains all differential operators,  their parametrices,  and  is invariant under composition (when defined) and formation of  adjoints.  The homogeneous  function  $\rho_s$ will be  called the {\it  principal symbol} of the operator.   It is uniquely defined by \eqref{symbol}. Moreover, the principal symbol  of a composed operator is the composition of the principal symbols.  Much as in the case of differential operators,  a pseudo-differential operator with symbol of the form  \eqref{symbol}   is elliptic if  and only if its principal symbol  $\rho_s(x,\xi)$ is invertible for $\xi\neq 0.$ Moreover,  the formula \eqref{par} for the parametrix extends to this class.

\smallskip

Let us discuss now the Agranovich reduction. 

The index bundle $\Ind (L,B)$ of a family of elliptic boundary value problems  coincides with the index bundle  of the family of operators defined  by the leading terms of operators  $\cLL(x,D)$ and $\cB_\la(x,D)$ respectively. Indeed, the  linear deformation of  lower order terms to $0$ produces  a homotopy between the corresponding Fredholm operators induced on Hardy-Sobolev spaces.  Therefore, with regard to  the computation of  $\Ind (L,B) $ we can safely assume that both $\mathcal L$ and $\mathcal B^1,\dots, \mathcal B^r$ are  homogeneous  polynomials  of  degree $k$ and $k_i$ respectively, which we will do from now on. In particular by $A_3$ we have that $\mathcal B_{\la}$ is independent of $\la.$

If $K$ is  the compact set  arising  in 
  assumption  $A_3,$  then for any $ x\in\Omega-K$  we have:
 \begin{equation}
\label{infty}
 \cLL(x,D) = \mathcal{L}_\nu(x,D). 
\end{equation} 
Being  ellipticity  an open condition, we can extend   $\cL $ to a   parametrized  family  of  elliptic  operators  (again denoted by $\cL$ ) defined on a open neighborhood $\co$ of $\bar \Omega $  and such that  \eqref{infty} still holds in $\co- K.$

For $u$ of compact support in $\Omega$ we have :\begin{equation}\label{pseudotop}
\cLL(x,D) u= (2\pi)^{-n} \int _{\R^n} e^{ix\xi} p(\la, x,\xi) \hat{u}(\xi) \, d\xi,
\end{equation}
 where $p$ is the principal symbol of the family $\cL.$

Let 
  \begin{equation}
\label{ops}
\tilde s(\la,x,\xi) =  \phi(|\xi| )p(\la,x,\xi) p(\nu,x,\xi)^{-1} + (1-\phi(|\xi|))\Id,
\end{equation}
where $\phi$ is  as in \eqref{par}.

By $A_3$, for $x \notin K,$  $p(\la,x,\xi) = p(\nu,x,\xi).$  Therefore, defining  
 $\tilde s(\la,x,\xi)=\Id$ outside of $\co$  we  can extend \eqref{ops} to a continuous map $\tilde s \colon \Lambda \times \R^{2n}\rightarrow GL(m;\C).$
 
  Each $\tilde s_\la$ is a symbol of order $0$ on $\R^n$  and, by the very definition of the topology in  $S^0(\R^n),$   $\tilde s$ is a continuous family of symbols such that $\tilde s_\la(x,\xi)= \Id $   for   $x \notin K.$

 Let  $\tilde \cSL$  be the operator associated  by  \eqref{pdo}  to  the symbol $\tilde s_\la.$ 
Then $\tilde \cS = \{\tilde \cSL \}_{\la\in \Lambda}$ is a  family  of pseudo-differential operators  on $\R^n.$ 

 It follows from \eqref{ops} that  the principal symbol of the family is given  by 
\begin{equation}
\label{syms}
\si(\la,x,\xi) = p(\la,x,\xi) p(\nu,x,\xi)^{-1}
\end{equation}    
for $x\in K$ and is the  identity   at points $(\la,x,\xi)$ with $x\notin K.$ Moreover, $\sigma$ extends in an obvious way to a  map defined on $\Lambda\times (\R^{2n}-K\times\{0\})$
with values in  $GL(m;\C).$

We will modify the family $\tilde\cS$ to a family of pseudo-differential operators with the same principal symbol but which has the property of being the "identity at infinity". For this, let $\psi \colon \Omega \rightarrow [0,1] $ be a smooth function which is identically 1 on $K$ and with compact support  $ K_1\subset \Omega $ and let 
\begin{equation}\label{psdo}
\cSL= \psi \tilde \cSL \psi + (1-\psi^2)\Id.
\end{equation}
 By the composition property,  the principal symbol of $\cSL$ is still  the same map $\si$ defined  
in \eqref{syms} and therefore each $\cSL$ is elliptic.  But now, being $\psi \equiv 0$ outside of $K_1,$ we have   

\begin{equation}\label{seeley}
 [\cSL u](x)= u(x)\, \hfor\, x\notin K_1. 
\end{equation}
Moreover, it is easy to see that the adjoint operator $\cSL^*$ has the same property.

The class of elliptic pseudo-differential operators such that both  the operator and its adjoint  verify  \eqref{seeley}   was introduced by Seeley in \cite {[Se]}. It plays a central role in the proof of the index theorem in \cite{[At-Si-1]}. We will denote this class of operators with $\pc.$  By \cite[Theorem 1,\, section\,1.2.3.5.] {[R-S]} each operator $\cq \in \pc $  extends to  a bounded operator $Q$ from $\hsn$ into itself. Moreover,  the correspondence sending the  symbol $\rho$ of the operator to the induced operator Q on $\hsn$ is  a continuous map from $S^0(\R^n)$ into $\mathcal{L}(\hsn)$ endowed with the operator norm topology.  
Taking into account our previous discussion,  the family $\cS $ defined by \eqref{psdo} induces a family of bounded  linear operators $ S\colon \Lambda \rightarrow\mathcal{L}(\hsn).$
If $(\cL,\cB)$ is a smooth family of boundary value problems parametrized by a smooth manifold $\Lambda$,  then the partials of the symbol of $\cS$ with respect to the coordinates $\la_i$ of $\la$ admit bounds of the form \eqref{boundpd}.  Therefore, $S$ is  a smooth family  whenever $(\cL,\cB)$ is  smooth.
\smallskip

  The following theorem is a version of the  Agranovich reduction \cite[Theorem 17.4]{[Ag]} for families of elliptic boundary value problems. 
  
  \begin{theorem}\label{th:a2} Let $(\cL,\cB)$ be a continuous  family of  boundary value problems verifying  
  assumptions $ A_1$ to $A_3,$ then the  family $ S\colon \Lambda \rightarrow\mathcal{L}(\hsn)$ defined above 
  is a family of Fredholm operators of index $0$  and 
  \begin{equation}
\label{ agranred}
\Ind (L,B) = \Ind S.
\end{equation}
 
\end{theorem}

\proof 
We will need to compare the operators  $\cS_\la$ and  $\cL_\la$.  The latter
is elliptic only on $\co$ and may not have  an elliptic extension to all of $\R^n.$ This problem can be handled by constructing a compact manifold to which $\cLL$ extends, but we prefer to avoid this  construction and instead  we choose a  compact neighborhood  $W$ of $\bar \Omega $ in $\co$ and we notice that, by \eqref{seeley},  $S_\la$ sends $\hsw$ into itself.  We will consider  $S_\la$ both as a bounded operator on $\hsn$ and on $\hsw$ and will split the proof of the theorem \ref{th:a2}  into a sequence of lemmas: 

\begin{lemma}\label{le:a3}  $S_\la\colon \hsw\rightarrow\hsw$  is Fredholm of index $0.$ Moreover, $S_\la L_\nu -L_\la  $ is a compact operator from $\hsw $ into itself. 
\end{lemma}
\proof
Each $S_\la$ and, as a matter of fact,  any  elliptic  pseudo-differential operator $\mathcal Q \in\pc$ has a parametrix $\mathcal P$  of  the same form.  Being $\Id -\mathcal{P} \mathcal{Q} $ of order $-1,$  the induced operator $\Id - PQ\colon \hsw\rightarrow \hsw$ factors through $H^{s+1}_W(\co,\C^m).$ Since $H^{s+1}_W(\co,\C^m) $ is compactly embedded in $\hsw,$   it follows that  $P Q $ is a compact perturbation  of the identity  and  moreover the same holds for $QP.$ Therefore, 
 $Q\colon\hsw\rightarrow \hsw$ is Fredholm by a classical  characterization of Fredholm operators. Since $S_\nu =\Id ,\ \ind S_\la=0\, \hforall \, \la .$ The second assertion follows again from the compact embedding of  $H^{s+1}_W(\co,\C^m)$ into $\hsw$ and the fact that the principal symbol 
  of $\cLL$ coincides with the principal symbol of $\cS_\la \circ \cL_\nu$ by the composition property. \qed\vskip5pt

\begin{lemma}\label{le:a2} Each  $S_\la\colon \hsn\rightarrow\hsn$  is Fredholm. Moreover, the index bundles of $S$ viewed  either as a family of Fredholm operators on  $\hsw$ or as a family on $\hsn$ are the same. 
\end{lemma} 
\proof 
We have a commutative diagram 
\[\bfig
\morphism (-500,0)<200,0>[0`; ]\morphism (-500,500)<200,0>[0`; ]
 \square|arra|<700,500>[\hsw`\hsn`\hsw`\hsn; i` S_\la`S_\la `i]
\square(700,0)|arra|/>``-->`>/<1000,500>[\hsn`\hsn/i\hsw`\hsn`\hsn/i\hsw;\pi```\pi]\morphism (2300,0)<200,0>[`0; ]
 \morphism (2300,500)<200,0>[`0; ]\efig\]
Being the support of $(S_\la u -u)$ contained in $W,$ the vertical dashed arrow induced by  $S_\la$ in the quotient spaces coincides with the identity.  Since an exact sequence of Hilbert spaces splits  into  a direct sum and since direct sums of Fredholm operators belong to the same class, each  $S_\la\colon\hsn\rightarrow\hsn$ is Fredholm.  Also the second assertion follows from the above diagram and the additivity of index bundle. \qed\vskip5pt

  Let  us take a bounded  extension operator  $E\colon \hs\rightarrow H^{s}_W (\co, \C^m)$  such that the values of $Eu$ on $\co -\bar \Omega$ depend only on the values of $u$ on $\bar \Omega-K_1.$  In order to obtain  such an operator $E$   it is enough to consider  the extension from $\hk$ to $\hkn$ constructed in \cite[ section 3.6 ]{[Ag]}, which verifies the above property, multiplied by a smooth function which coincides with $1$  on  $\bar\Omega$ and  with  support in $W.$ Finally, let  $R\colon\hsw\rightarrow\hs$ be  the restriction operator.

\begin{lemma}\label{le:a1} 
 $R S_\la E  L_\nu -   L_\la \colon  \hk\rightarrow \hs $ is compact  for all $\la \in \Lambda. $
\end{lemma}  
\proof
Here we closely follow the arguments used in the proof of   \cite[ Theorem 17.4]{[Ag]}. Since we are working with a different class of operators, we  include the proof  for the sake of completeness. 

 We will first show that:
\begin{equation}
\label{eq:era}
  S_\la(ER -\Id) = ER-\Id \hand   (ER -\Id) S_\la= ER-\Id.
\end{equation}

Indeed, denoting with $ \tilde  S_\la$ the operator induced by $\tilde \cS_\la$ on Hardy-Sobolev spaces, by definition of $\cS_\la,$  we have: 
\[ S_\la E  R -S_\la=  \psi \tilde S_\la \psi (ER-\Id ) +  (1-\psi^2) (ER-\Id).\]
 
 But $ \psi \tilde S_\la \psi (ER-\Id)=0$ because the support of the function $\psi $ is contained in $\Omega$  and $ (1-\psi^2) (ER-\Id)= ER-\Id $ by the same reason. This proves the first relation in \eqref{eq:era}. The proof of the second relation  is similar. 

Applying R to the first equation in \eqref{eq:era} we get 
 \begin{equation}\label{rser}
R S_\la E  R=  R S_\la \, \hforall \, \la. 
\end{equation}
Let us represent  $L_\la$ defined on  $\hk$ in  the form $L_\la= RL_\la E,$  where the $L_\la$ on the right hand side is viewed as an operator on $\hsw.$  Using \eqref{rser}  \begin{equation}\label{00}\begin{array}{l}
 R S_\la E   L_\nu -   L_\la=RS_\la ERL_\nu E -  RL_\la E\\=RS_\la L_\nu E - RL_\la E=R(S_\la L_\nu - L_\la)E 
\end{array}
\end{equation}
is compact by the second assertion in lemma \ref{le:a3}. \qed\vskip5pt

\begin{lemma}\label{le:a4}
The family  $  S'$ defined by $S'_\la=RS_\la E\colon \hs\rightarrow \hs $ is a family of Fredholm operators  and $ \Ind  S' =  \Ind S.$
   \end{lemma}

\proof We will first show  that $S'_\la $ is Fredholm. 
  Using \eqref{seeley} if $S'_\la u_n\rightarrow f,$ then  $u_n$ restricted to $ \bar \Omega-K_1$ converges in  $H^s(\Omega -K_1)$ to the restriction of $f.$  By the construction of $E,$    $Eu_n\rightarrow Ef$ in $H^s_W(\co -K_1).$   It follows that  $S_\la Eu_n\rightarrow Ef.$  Since $S_\la$ has a closed image, there exist a $w\in\hsw$ such that  $S_\la w= Ef.$ But then the restriction of $w$ to $\co -K_1$ coincides with $Ef$  which  implies that  $ERw=w$ and hence  $RS_\la ERw = RS_\la w=f.$ This shows  that $\im S'_\la$ is closed.
 
 Applying $R$ to the left of the first equation in \eqref{eq:era}  and $E$ to the right of the second we get
\begin{equation}
\label{eq:comm}
  S'_\la R = RS_\la  \hand   ES'_\la = S_\la E.
\end{equation}

 The second equation shows that $E$ sends $\ker S'_\la$ into $\ker S_\la$ and since $E$ is injective, $ \dim \ker S'_\la$ is finite.  In order to show that $\dim \coker S'_\la$ is finite we observe that  the first equation in \eqref{eq:comm} shows that $ R\colon \hsw\rightarrow \hs$ sends $\im S_\la$ into $\im S'_\la$ and hence induces $\bar R \colon \hsw/ \im S_\la\rightarrow \hs/\im S'_\la.$  Being $R$ surjective, the same holds for $\bar R$ and therefore $\dim \coker S'_\la$ is finite. 

 Let us show now that $\Ind S'=\Ind S.$ If $F$ is a finite dimensional subspace of $\hsw$ such that 
$\im S_\la + F = \hsw \hforall  \la \in \Lambda,$ then $H= ER(F)$ enjoys the same property because $(ER-\Id)(F) \subset \im S$ by \eqref{eq:era}.  

Applying $R$ to both sides  we get
 $$\im RS_\la +R (H) = \hs \hforall  \la \in \Lambda.$$ But, by the first equation in \eqref{eq:comm},
 $\im RS_\la \subset \im S'_\la,$ which shows that $H'= R(H)=R(F)$  is transverse to $\im S'_\la \hforall \la.$  Notice also that $E$  sends  isomorphically  $H'$ into  $H$ with inverse $R.$   Denoting with  $G_\la$  and  $G'_\la$  the inverse images of $H$ and $H'$ under $S_\la$ and  $S'_\la$ respectively, the second equation  in \eqref{eq:comm} implies that  $E(G'_\la ) \subset G_\la.$  On the other hand, being $E$ injective and since   $$\dim G_\la =\dim H = \dim H' = \dim G'_\la,$$ it follows that $E$ induces a vector bundle isomorphism between vector  bundles $G'$ and $G$ over $\Lambda.$ Thus 
$$ \Ind S' = [G'] -\theta(H') = [G] - \theta(H)=  \Ind S.$$\qed\vskip5pt

Now we can complete the proof of theorem \ref{th:a2}.  Let $(\bar L_\la, \bar B_\la)_{\la\in \Lambda}$ be the family of operators  defined  as   the composition 
 \[\bfig
 \morphism(0,0)<1100,0>[\hk`\hsb;(L_\nu,B)]
 \morphism(1100,0)<1400,0>[\hsb`\hsb;S'_\la\times\Id]
 \efig.\]
By logarithmic property of the index bundle,  $$\Ind (\bar L,\bar B) = \Ind (L_\nu,B) + \Ind (S'\times\Id) =  \Ind (S'\times\Id),$$ 
since both $L_\nu$ and $B$ are independent from $\la.$ On the other hand, by lemma \ref{le:a1},  $L -\bar L $ is a family of compact operators. Hence,  so is  $(L,B) - (\bar L,\bar B),$ and therefore
$$ \Ind(L,B) = \Ind (\bar L,\bar B)= \Ind (S'\times\Id)= \Ind S' =\Ind S,$$  by lemmas \ref{le:a1} and  \ref{le:a4}.\qed\vskip5pt

\vskip20pt
\newpage

\subsection {Proof of the bifurcation theorems \ref{th:40} and \ref{th:50}.} 
 \subsubsection{Proof of theorem {\rm \ref{th:40}.}}
    It follows from  the discussion  in the second part of  Appendix B that, for $s>n/2,$ the family of  nonlinear differential operators  \[(\cf,\cg) \colon \R^q \times C^\infty( \Omega;\R^m) \rightarrow C^\infty( \Omega;\R^m) \times\ C^\infty(\partial\Omega;\R^r)\]
  induces a smooth map 
\begin{equation}
\label{nbvp1}
 h=(f,g) \colon \R^q \times H^{k+s}(\Omega;\R^m ) \rightarrow    \hsbr,
\end{equation}
with $\R^{q} \times\{0\}$ as a trivial branch. 
The Frechet derivative  of $h_\la$ at $ 0$   is the operator 
$(L_\la,B_\la)$ induced by the linearization $(\cL_{\la}, \cB_\la)$  at  $u\equiv 0.$ Since, for any $\la\in R^q,$ $(\cL_{\la}, \cB_\la)$ is elliptic, using proposition \ref{prop:ell1}, we can find a neighborhood  $O$  of $0$ in $\hk$  such that  $h\colon \R^{q} \times O \rightarrow  \hsbr $   is a smooth  family of  semi-Fredholm maps. 

 By  $H_2,$  the family  of boundary value problems $(\cL, \cB)$ extends to a smooth family  parametrized by $S^{q}$  which clearly verifies the assumptions $A_1$ to $A_3$ of  section \ref{sec:4.1}  with $\nu = \infty \in S^{q}.$  Hence, the induced family on Hardy-Sobolev spaces also  extends  to a smooth family $(L,B) \colon S^q \ra \mathcal L(\hkr, \hsbr)$  Moreover, $(L_{\infty},B_{\infty})$ is  invertible by $A_2.$ Thus  $(L_{\la},B_{\la})$  is Fredholm of index $0,$ for all $\la\in S^{q}$  and, by continuity of the index of semi-Fredholm operators,  the map  $h\colon \R^{p} \times O \rightarrow  \hsbr $  is a smoothly parametrized family of Fredholm maps of index $0.$
 
 In order to simplify our notations, in the rest of this section  we will abbreviate $(L,B)$ to  $L$ when no confusion  arises. 

 Since $h$ is defined only on the open subset $\R^q$ of $S^q$  we cannot apply directly theorem \ref{th:1.1} to $h$  in order to find a bifurcation point.  Instead we will  use the assumption  $H_2$  in order to compute the local  index  $\beta(h,\R^q)$ from the family index theorem applied to  $L=(L,B) \colon S^q \ra \mathcal L(\hkr, \hsbr).$

Since $L_\la$ is invertible in a neighborhood of  $\infty\in S^q, $ the pair $(h,\R^q)$ is admissible and  the local bifurcation index $\beta(h,\R^q) \equiv \beta_{S^q}(h,\R^q)$ is defined. Being  $L$ an extension of  $Dh_{-}(0)$ to all of $S^q,$  by the very definition of the local bifurcation index,  $\beta(h,\R^q)= \J(\Ind L).$  

If, under the hypothesis of theorem  \ref {th:40},  we can show that $\J(\Ind L) \neq 0 $ in $\J(S^q),$  then  the family $h$ must have a bifurcation point  $\la \in \R^q,$   by  $\text{B}_1.$ This would  complete the proof of the theorem,  since by proposition \ref{prop:boots}   it follows  that a bifurcation point of the map  $h$  is also a bifurcation point for smooth classical solutions of  \eqref{bvp1} in the sense of definition \ref{def:1}.
\vskip 10pt

The remaining part of the proof is devoted to show that $\J(\Ind L) \neq 0 $ in $\J(S^q).$  For this, we are going to to compute $\J (\Ind L)$ from  the degree of $\sigma$ using  the complexification $L^c$ of $L.$   
Since  $\ker L^c =\ker L\otimes \C ,$  from definition of the index bundle in \eqref{defind}  it follows that 
\be\label{compl} \Ind L^c= c(\Ind L),\ee  where  $c\colon \tilde KO \rightarrow \tilde K$ is the complexification homomorphism.
 
 By Bott periodicity, $\tilde K(S^q)=0 $ for $q$ odd, while for  $q=2 k,$  $\tilde K(S^q) $ is an infinite cyclic group. It is generated by $\xi_q =\left( [P^1(\C)\times\C ]-[H]\right)^{k},$  where $H$ is the tautological  line bundle over the complex projective space $P^1(\C)\cong S^2.$  On the other hand,  the periodicity theorem for $\tilde{KO}$ gives  $\tilde{KO}(S^q) \cong  \Z$  for $q\equiv 0,4 \mod 8,$   $\Z_2$ for  $q \equiv 1,2 \mod 8$  and vanishing in all remaining cases.  
  From the homotopy  sequence of  the fibration of classifying spaces for  $\tilde{KO}$ and $ \tilde{K}$ (see \cite[ section 13.94]{[Sw]}) 
it follows  that  $c\colon\tilde{KO}(S^q) \rightarrow \tilde{K}(S^q)$ is an isomorphism for $q\equiv 0 \mod 8 $  and a monomorphism with image generated by $2\xi_q$ for  $q\equiv 4 \mod 8.$  

For $q=4s,$ we take as generator of $\tilde{KO}(S^q)$ an element $\nu_q$ such that
 \be\label{gen} c(\nu_q)= \begin{cases} \xi_q  &\hif q\equiv 0 \mod 8 \\ 2\xi_q  &\hif  q\equiv 4 \mod 8.\end{cases}\ee 
 
  With this choice of generators, each  element $\eta \in\tilde K(S^q)$ with $q=2k$ is  uniquely determined by its {\it degree}  $d(\eta) \in \Z$  verifying  $ \eta = d(\eta)\, \xi_q, $ and, for $q =4s,$ each  element  $ \eta$ of $\tilde{KO}(S^q)$  has  a degree defined in the same way. 
 
 By \eqref{gen}, for any $\eta \in  \tilde{KO}(S^q),$   
 \be\label{degr} d(c(\eta)) = \begin{cases} d(\eta) &\hif q\equiv 0 \mod 8 \\ 2d(\eta)  &\hif  q\equiv 4 \mod 8.\end{cases}\ee
 
The degree of an element  $\eta \in \tilde{K}(S^q)$ can be computed as a characteristic number in several ways.   We will use the Chern character $ \ch \colon \tilde{K}(-) \rightarrow  \h^{*}(-;\C)$  with values in de Rham cohomology with coefficients in $\C$  which is adequate to our purposes.   By Bott's integrality theorem,  $ ch=ch_k\colon \tilde{K}(S^{2k}) \rightarrow  \h^{2k}(S^{2k};\C)$ is injective with image given by $\im \ch_k = \Z\, u_{2k},$  where  $u_{2k} = \ch_k (\xi_{2k})$ is the class  of  the volume form of $S^{2k}$ \cite[Chap.18, Theorem 9.6] {[Hu]}.  Hence, for any   $\eta \in \tilde{K}(S^{2k}),$ we have  $ \ch( \eta ) = d(\eta)u_{2k}$ and therefore 
\be\label{cher} d(\eta) = < \ch( \eta); [S^{2k}]> ,\ee  the right hand side being  the evaluation of  $\ch(\eta) $ on the fundamental class $[S^{2k}]$  of the sphere.

Since  the complexified family  $(\cL^c,\cB^c )$  verifies  assumptions $A_1$ to $A_3$ of theorem \ref{th:a2},  $\Ind L^c=\Ind S,$ 
 where $S $ is  induced by  the family  of pseudo-differential operators $\cSL$ defined by \eqr{psdo}.  By Fedosov's formula (see Appendix C)  with $j= n+k,$ we get
 
 \be \label{chern} \ch( \Ind S )= K_k   \oint _{S^{2n-1}} tr(\sigma^{-1}d\sigma)^{2(n+k)-1},\ee  where $\oint$ denotes the integral along the fiber and $K_k =-\displaystyle{\frac{( n+k-1)!} {(2\pi i)^{n+k}(2n +2k-1)!}}.$ 
 \medskip
 
The evaluation on the fundamental class in de Rham cohomology  is given by integration over the sphere. Hence, using Fubini's theorem for integration along the fiber,  from \eqref{chern} we get
\begin{equation}\label{eqfin}
 d(\Ind L^c) = K_k \int_{S^{2k}}  \oint _{S^{2n-1}} tr(\sigma^{-1}d\sigma)^{2(n+k)-1} =  K_k \int_{S^{2k}\times S^{2n-1}} tr (\sigma^{-1}d\sigma)^{2(n+k)-1},
\end{equation} 
where  the right hand side is the ordinary integration of  the $(2k+ 2n-1)$-form $tr (\sigma^{-1}d\sigma)^{2(n+k)-1}$ over $S^{2k}\times S^{2n-1}.$ 

 Thus $d(\Ind L^c)$ coincides with  $d(\sigma)$ defined in \eqref{fed}. Using  \eqr{degr} we obtain
 \be\label{degre} d(\Ind L) = \begin{cases} d(\sigma) &\hif q\equiv 0 \mod 8 \\ \frac{1}{2}d(\sigma)  &\hif  q\equiv 4 \mod 8.\end{cases}\ee
 On the other hand, for  $q =4s,$  $J(S^q)\simeq  Z_{m(q/2)}$ and  $J(\Ind L) =0$ if and only if  $d(\Ind L)$ is divisible by $m(q/2).$ Now, theorem \ref{th:40} follows from \eqref{degre} and the definition of  $n(q)$ in \eqref{muca}.\qed\vskip5pt

\vskip20pt

\subsubsection{Proof of theorem {\rm  \ref{th:50}}} 
Let us first recall the clutching construction. 
 Given a continuous  map $ G \colon S^{q-1} \rightarrow GL(m;\C),$  taking  two trivial complex bundles of rank $m$  over the upper and lower  hemispheres $D_\pm$ of $S^q$   we obtain  a bundle $\eta_G $ over $S^q$ by identifying  $(\la,v)\in  \partial D_+\times \C^m$ with  $(\la,G_\la v)\in  \partial D_-\times \C^m.$   The  isomorphism class of  $\eta_G$ depends only on the homotopy class of $G.$ Moreover, the clutching construction  extends to an isomorphism 
between  $\pi_{q-1}{GL(\infty;\C)}$ and  $\tilde K(S^q)$. An analogous construction establishes an isomorphism  of  $  \pi_{q-1}{GL(\infty;\R)}$ with  $\tilde{KO}(S^q)$  which coincides with the inverse of the isomorphism $\partial_0$ in lemma \ref{lem:bell}. 

  If the map $G$ is smooth, choosing  appropriate connection-forms  on $D_\pm,$ one can compute  $d(\eta_G)$ as

 \begin{equation}\label{degeta} d(\eta_G)=<ch_k(\eta_G);[S^{2k}]> =  \displaystyle {\frac{-(k-1)!} {(2\pi i)^{k}(2k-1)!}} \int_{ S^{2k-1}} tr (G^{-1}dG)^{2k-1}.\end{equation} 
 A proof of this can be found  in  section 3.2 of \cite{[Fe-2]} (see also \cite{[Al-Yo]} in the the real case).
 
For $q=4s,$  let  $\la_0=0$  be an isolated singular point of $L.$  Without loss of generality  we can assume that ${\la} _0$  is the north pole of $S^q$ and that the open neighborhood  $U$ isolating $\la_0$ from the rest of $\Sigma(\bar L)$ contains  the upper  hemisphere $D_{+}.$ 

We extend $L \mid_{D_{+}}$  to  a  family $\tilde L$ defined on all of $S^q$ such that  
$\tilde  L_{\la}$ is an isomorphism for $\la\in D_{-} .$  If $A_{+}$ is  any   parametrix for $ L_{+}$  and if  we  take as $A_{ -}= \tilde  L_{-}^{-1},$  then, arguing as in the proof of  proposition \ref{pr:alex},  we can show that the homomorphism $\partial $ of the diagram \eqref{eq:gf} sends $\tilde L$  to the  family of matrices $N$ whose stable homotopy class is taken as definition of  $\ga_f$ in  section \ref{sec:3.3}.  By  commutativity of the diagram  \eqref{eq:gf}  and since the clutching construction is the inverse of $\partial_0,$ we have 
\begin{equation}  \label{finaly} \Ind( L,U)=  \Ind \tilde L= [\eta_N]. \end{equation}  
As in the proof of theorem \ref{th:40}  we can compute $d(\Ind( L,U))$ from the complexification of $[\eta_N].$  It is easy to see that   $c[\eta_N]$  is the vector bundle associated by the clutching construction to the complexification $N^c$ of $N.$  By  \eqref{degr}, 
 \be\label{degIND} d(\Ind( L, U)) = \begin{cases} d(\eta_{N^c}) &\hif q\equiv 0 \mod 8 \\ \frac{1}{2}d(\eta_{N^c})  &\hif  q\equiv 4 \mod 8.\end{cases}\ee
   
By proposition \ref{prop:ls},  $R$ is homotopic to $N$  and hence from \eqref{degeta} we obtain 
\begin{equation}\label{degR} d(\eta_{N^c})=d(\eta_{R^c})=  \displaystyle (-1)^{s+1}  {\frac{(2s-1)!} {(2\pi)^{2s}(4s-1)!}} \int_{ S^{4s-1}} tr ({R^c}^{-1}dR^c)^{4s-1}.\end{equation}
The right hand side of \eqref{degR} coincides with the degree $d(\la_0)$ defined in \eqref{fedlo} because $ tr ({R^c}^{-1}dR^c)^{4s-1}= tr ({R}^{-1}dR)^{4s-1}.$ This gives 
\be\label{degree}d( \Ind( L,U))= \begin{cases} d(\la_0)&\hif q\equiv 0 \mod 8 \\ \frac{1}{2}d(\la_0) &\hif  q\equiv 4 \mod 8.\end{cases}\ee

Now, the  assertion  $i)$ follows from  \eqref{degre},  \eqref{degree} and the additivity property  \eqref{eq:2.7} of the index bundle.  Under the isomorphism $\J(S^{q})\simeq \Z_{m(q/2)},$  $ \J(\Ind(L,U)) $ coincides with mod $m(q/2)$ reduction of  $d\left(\Ind( L,U)\right).$  Thus the first part of ii)  follows from  \eqref{degree}, the definition of $n(q)$  and  $\text{B}_{1}.$   For the second part it is enough to observe that if  $d(\sigma)-d(\la_{0})$ is not a multiple of $n(q),$ then  $\beta(h,\Lambda -\{\la_{0}\})\neq 0 $ in  $J(S^{q}),$ by  additivity of the bifurcation index.

 \vskip 20pt
    
 \vskip20pt
\newpage
 
\specialsection{ \footnotesize APPENDIX}\label{sec:5}

\subsection{A. Properties of the index bundle}  \label{sec:5.1}

Since our construction  of the index  bundle  differs from the one in \cite{[At],[Ja]},
we briefly describe the proofs of its properties.

\begin{proposition} 

The index bundle $\Ind L$ verifies:  
\begin{itemize}

\item[i)]\emph{ Functoriality:}  If $L\colon \Lambda\rightarrow  \Phi (X,Y)$ be a family of Fredholm operators  and  
$\alpha \colon \Sigma \rightarrow  \Lambda$ is a continuous map between compact spaces, then 
 \[\Ind L\circ\alpha = {\alpha}^*(\Ind L),\]
where $\alpha^*\colon KO(\Lambda )\rightarrow KO (\Sigma)$ is the
 homomorphism induced by $\alpha.$
 
\item[ii)] \emph{Homotopy invariance:}  Let $H\colon[0,1]\times   \Lambda\rightarrow  \Phi (X,Y)$ be a
homotopy,   then $\Ind  H_0 =\Ind  H_1.$

\item[iii)] \emph{Additivity:} $ \Ind  \bigl( L\oplus M \bigr) = \Ind L + \Ind M.$

\item[iv)] \emph{Logarithmic property:}  $\Ind\bigl( LM \bigr) = \Ind L + \Ind M.$

\item[v)] \emph{Normalization:} If  $L$ is homotopic to a family in $GL(X,Y),$  then $\Ind L =0.$ Moreover,  the converse holds if  $Y$ is a Kuiper space.
 \end{itemize}
\end{proposition}

\proof Taking the same subspace $V$ in the definition of the index bundle for both  $L$ and $L\circ \alpha$,  property  i) follows plainly from  the  definition of $\alpha^*(E)$.  Now, ii) follows from  i) applied to the top and bottom inclusions of $ \Lambda$  in $[0,1]\times  \Lambda.$ The proof of iii) is straightforward. Assuming $X=Y=Z,$  iv)  reduces to iii) thanks to a well known homotopy  between $ Id\oplus LM$ and  $L\oplus M$ \cite[Theorem $7.2$]{[Bo-Bl]}. The general case follows easily from this.  Another way to prove iv)   is by observing   that in the construction of the index bundle one can take instead  of a  finite dimensional subspace   $V$ of $Y$  any  finite dimensional subbundle of  $\Lambda \times Y$ transverse to $L.$   Now, if $\Theta(V)$ is transverse to $LM,$ then  $\Theta(V)$ is transverse to $L$ and $E = L^{-1}\Theta(V) $ is transverse to $M$. Then, denoting by $ F =  M^{-1}E$, in $KO(\Lambda )$ we have
\begin{equation}
 \Ind  \bigl( L M \bigr) = [F]- [\Theta(V)] =  ([F]- [E]) + ([E]- [\Theta(V)] ) = \Ind  L + \Ind M.
\end{equation}
 The proof of  v) can be found in \cite[Theorem $1.6.3$]{[Fi-Pe]}. 
\vskip20pt

\subsection{B. Elliptic boundary value problems}  \label{sec:5.2}
\subsubsection{Linear elliptic boundary value problems}  \label{sec:5.2.1}
We begin with a brief summary of the relevant linear  theory. We will  work over the field $\C$ of complex numbers considering   real coefficients as a special case.  For nonlinear systems it becomes  natural  to take the opposite viewpoint. 

For $\alpha  = (\alpha_1,\ldots 
,\alpha_{n})$  an $n$-tuple of nonnegative integers, we set
\[
  D_j = i^{-1} {\frac{ \partial }{ \partial x_j}},\,  D^{\alpha }=\prod _{i=1}^{n}\left( D_i\right)^{\alpha_ i}, \vert \alpha \vert  
=\sum^{n}_{i=1}\alpha_i \, \hbox{and for }\, \xi  \in  \C^{n}, \, \xi^{\alpha }=\prod_{i=1}^n {\xi_i}^{\alpha_i}.
\]

 Let $\Omega$ be an open bounded subset of $\R^n$ with smooth boundary.  We will consider partial differential operators acting on smooth vector functions $u\colon \Omega\rightarrow \C^m$ of the form

\begin{equation}\label{ellop}
 {\mathcal L}(x,D)u=\sum^{}_{\vert \alpha \vert \le k}a_{\alpha }(x)D^{\alpha }u(x),
\end{equation}  
where   $a_{\alpha }\in C^\infty (\bar\Omega;  {\C}^{m\times m}).$ 
 The {\it principal part} of $\mathcal L$ is the expression \eqref{ellop} containing only the leading terms with $|\alpha| =k.$ 
The {\it  principal symbol} of $\mathcal L $ is the matrix function  $p$  defined on  $\Omega\times \C^n$ by
\begin{equation}\label{ symb} p(x,\xi ) \equiv \sum^{}_{\vert \alpha \vert =k}\xi ^{\alpha }a_{\alpha }(x). \end{equation}

The operator $ \cL(x,D)$  is called elliptic  if its principal symbol
 verifies 
\begin{equation}
\label{ellipt} \det \ p (x,\xi )\neq 0   \, \hbox{ for all }x \in\bar\Omega,\,  \xi \in \R^{n} -\{ 0\}.\end{equation}

 $ \cL(x,D)$ is called {\em properly elliptic}  if $km=2r$ and  for any  $x\in \partial \Omega$ and any vector $\xi \neq 0 $ tangent to the boundary at $x,$ denoting with $\eta$  be the  inward normal to $\partial \Omega$ at $x,$ we have that the polynomial  $\mathrm{det}\, p (x,\xi + z \eta )$ has exactly $r$ roots in the upper half-plane $ \Im z > 0.$  If we introduce  coordinates $(y_1,\dots , y_{n-1}, y_n)$  at  $x$  such that  $\partial \Omega $ is defined in a neighborhood of $x$  by $y_n=0,$ then, in terms of the ordinary differential operator   $ p(y_1, \dots , y_{n-1},0, \xi,  i^{-1}{ \frac {d}{dt}} ),$ the above condition means that the  subspaces  $ M^\pm (x,\xi )$ of  $ L^2( \mathbb R_\pm ;\C^m)$  whose elements are exponentially decaying  solutions of the system $ p(y_1,\dots , y_{n-1},0, \xi,  i^{-1}{ \frac {d}{dt}})v(t)=0$ at $+\infty$  and $-\infty$ have dimension $r.$
 
Let  ${\mathcal L(x,D)}$ be  an  elliptic operator  of order $k,$ properly elliptic at the boundary and let  $ k_i , 1\, \le \, i\, \le \, r,$ be integers such that  $0 \le  k_i \leq k-1.$ 
We will consider  $r$ operators $\{\mathcal B^1(x,D),..., \mathcal B^r(x,D)\}$ of order $k_i.$
\begin{equation}\label{bop}
 {\mathcal B}^i(x,D)u=\sum_{\vert \alpha \vert \le k_i}b^i_{\alpha }(x)D^{\alpha }u(x), 
\end{equation} 
where   $b^i_{\alpha }\in C^\infty (\bar\Omega; {\C}^{1\times m}).$

 The {\it boundary operator} is the operator matrix $\mathcal B(x,D)$ whose  $i$-th row is  $\cB^i(x,D).$ 
Thus $\cB(x,D)=[\mathcal B^1(x,D),\dots,\mathcal B^r(x,D)]^t.$ 

 The principal symbol of the boundary operator $\mathcal B(x,D)$ is by definition the matrix function $p_b(x,\xi) $  whose $ i$-th row  is 
\begin{equation} \label {boundsymb}
p_b^{i }(x,\xi)  =  \sum^{}_{\vert \alpha \vert =k_i}\xi^{\alpha} b^{i}_{\alpha }(x).\end{equation}

The boundary  operator  $\cB$ verifies the  {\em Shapiro-Lopatinskij  condition} with respect to ${\cL}(x,D)$ if,  for each $x \in  \partial \Omega$   and $\xi  \in  \R^{n}\backslash  \{0\}$ belonging to $T_x\partial \Omega,$ the subspace $ M^+ (x,\xi )$  is isomorphic to $\C^r$ via the map  $ u  \mapsto [ p_b(y_1,\dots , y_{n-1},0, \xi,  i^{-1}{ \frac {d}{dt}})v](0).$

Since the condition involves only ordinary differential equations with  constant coefficients, it can be reformulated in purely algebraic terms \cite{[Ag-Du-Ni]} but we will not   use this formulation here. 

 \begin{definition}\label{EBVP} Given an open bounded subset $\Omega$ of $\R^n$ with smooth boundary, an  elliptic boundary value problem   on $\Omega$ is  a pair    $({\cL}, {\cB})$  where ${\cL} ={\cL}(x,D)$ is an  elliptic operator on $\Omega$,  properly elliptic at the boundary,   and the  boundary operator ${\cB}= \cB(x,D)$  verifies the Shapiro-Lopatinskij condition with respect to $\cL.$ \end{definition}
 
For  any manifold $M,$ with or without boundary, there is an  associated scale of Hardy-Sobolev spaces $H^{s}(M, \C^m ),  s \in \R$ \cite{[Ch-Pi]}. Every  $u \in  H^{s}(M ;\C^m)$  has a well defined restriction to $\partial M $ belonging  to  $H^{s-1/2}(\partial M;\C^m ) $ and continuously  depending on $u$. 
 When  $ s \in \N $ and $M=\Omega$ an open subset of $\R^n$ with smooth boundary, denoting with $D^{\alpha }u$  the distributional derivative, we have  
\[H^{s}(\Omega; \C^m ) = \{u \in  L^{2}(\Omega;\C^m )\vert D^{\alpha }u \in   L^{2}(\Omega; \C^m )\hforall \vert \alpha \vert  \le  s\}\]   with the norm $\Vert  u  \Vert_s = \sum_{|\alpha|\leq s} | D^\alpha u |_2.$

Let   $ \tau \colon C^{\infty} (\bar \Omega)\rightarrow C^\infty (\partial \Omega)$ 
 be the trace operator.  The operator    \[(\cL,\tau \cB)\colon  C^\infty(\bar \Omega; \C^m ) \rightarrow C^\infty(\bar \Omega; \C^m ) \times C^\infty (\partial \Omega; \C^r)\] extends  to  a bounded operator  
\begin{equation}\label{fredo}
( L,B) \colon H^{k+s}(\Omega; \C^m ) \rightarrow  H^{s}(\Omega; \C^m ) \times\hbb, 
\end{equation}      
where $\hbb$ denotes  $ \prod_{i=1}^r H^{k+s-k_i-1/2}(\partial \Omega; \C).$

For any  elliptic  boundary value problem  $({\mathcal L},\mathcal B)$ the following Schauder type estimate holds  \cite{[Ag-Du-Ni]}:

 there exists a constant $c > 0$  such that for any $  u\, \in \, H^{k+s}(\Omega )$  
\begin{equation}\label{schauder} 
\Vert u\Vert_{k+s} \le  c \left( \Vert L  (u)\Vert_{s} + \sum^{r}_{i=1}\Vert B_i(u)\Vert_{k+s-k_i-1/2} + \Vert u \Vert
_{s}\right).\end{equation}

It follows easily from the above estimate that the operator $(L,B)$  has finite-dimensional kernel and closed image. Namely,  $(L,B)$  is left semi-Fredholm.
\vskip20pt
\subsubsection{Nonlinear elliptic boundary value problems}\label{sec:5.2.2}  Denoting  with  $k^*$  the number of multi-indices $\alpha $ with $\vert \alpha \vert  \le k,$ the $k-$ jet extension 
$ j_{k} : C^{\infty}(\bar \Omega, \R^m ) \rightarrow  C^\infty(\bar \Omega ,\R^{mk^*})$  is defined by 
$$
(j_{k}u(x))_{ \alpha } = D^{\alpha }u(x)\hbox{ for}  \vert 
\alpha \vert  \le k.$$
Given  a continuous family of smooth maps 
$\cf : \Lambda \times\bar \Omega  \times  {\R}^{mk^*} \rightarrow  {\R^m}$  we will  informally write $ \cf(\la, x,u(x),\ldots ,D^{k}u(x))\,  \hfor\,  \cf(\la, x, j_{k}u(x)).$
 As in the case of linear differential operators we will not distinguish in the notation the map $\cf$ from the  family of nonlinear  operators 
 $\cf \colon \Lambda\times C^\infty( \bar\Omega; \R^m)\ra  C^\infty( \bar\Omega; \R^m)$ defined by the above expression.
 However, we will use roman alphabet  to denote the corresponding operators induced in Hardy-Sobolev spaces.
 
 An argument based  Sobolev's embedding theorems, shows that for $ s > n/2$ (which we will always assume)  the family  $ \cf(\la, x,u(x),\ldots  ,D^{k}u(x))$  extends to  a  continuous family of  smooth maps  $f\colon \Lambda\times H^{k+s}(\Omega;\R^m )\rightarrow H^{s}(\Omega;\R^m ).$ Moreover, if $\Lambda$ is a smooth manifold and $\cf$ is smooth, so is $f.$ 
 
  Indeed, for any $\la\in \LL,$  the nonlinear operator $f_\la $ is the composition  of $j_k$ with  the Nemytskij operator associated to the map $\cf_\la.$  The operator $j_k$ extends to a bounded linear operator from $\hkr$ to  $\hsr$ while, for $ s > n/2,$   the associated Nemytskij  operator induces a smooth map from  $\hsr$ into itself (see \cite[ Theorem 11.3]{[Pa]}).    Moreover, the argument used in the proof of  \cite[Theorem 11.3]{[Pa]}  automatically gives the continuous dependence on parameters of the derivatives of $f_\la,$ if $\cf$ is a continuous family of smooth maps. The same argument allows to show that $f$ is smooth if so is $\cf.$

  Together  with  $\cf,$ we will consider $r$ nonlinear  boundary conditions  of  order  $k_i$  with $0 \le  k_i \leq k-1.$ These are defined by  $r$ continuous families of  smooth maps  $ \cg_i \colon \Lambda\times \bar\Omega \times  {\R}^{m{k_i}^*} \rightarrow  {\R }, \,1 \le  i \leq \hbox{r}.$
  
   Composing the obvious  projections from $\R^{mk^*}$ into $\R^{mk_{i}^*}$  with the functions  $ \cg_i,\, 1\leq i\leq r,$  we obtain a map  $\cg=(\cg_1\dots \cg_r)\colon\Lambda\times\bar\Omega \times \R^{mk^*} \ra \R^r$   and hence a family of nonlinear {\it boundary operators} $ \cg\colon \Lambda\times C^\infty(\bar \Omega;\R^m)\ra C^\infty(\partial \Omega;\R^r) $  defined  by   
 \begin{equation} \cg(\la, u)=  
\left( \ \tau \cg_1(\la, x,u,\ldots,D^{k_1}u),\dots,  \tau \cg_r(\la, x,u,\ldots  ,D^{k_r}u)\right),\end{equation}
where $\tau$ is the restriction to the boundary.

 The above discussion, together with the  well known continuity property  of the trace $\tau,$    allows to conclude that the map $(\cf,\cg)$ extends to  a  continuously parametrized  family of smooth maps  
\begin{equation}
\label{bvpa}
(f,g)\colon  \Lambda\times H^{k+s}(\Omega;\R^m ) \ra \hsbr .\end{equation} 

 For each  fixed $\la,$ {\it the  linearization} of  $(\cf_\la,\cg_\la)$ at a smooth function $w$  is the linear operator:
\begin{equation}
\label{lin}
\begin{array}{ c}
 {\mathcal L}_\la(x,D) u(x) =\sum^{}_{\alpha}a_{\alpha }(\la, x) D^{\alpha }u(x)\\
\mathcal {B}_\la(x,D) u(x) =\tau \sum_{ \alpha} b_\alpha (\la,x)D^{\alpha }u(x),
\end{array}
\end{equation}
\noindent where,  denoting by  $v_{j\alpha}$ the variable corresponding to $D^\alpha u_j,$  the  $ij$- entries of the matrices $a_{\alpha } \in C^\infty( \Lambda\times \bar\Omega; \R^{m\times m} )$ and  $b_{\alpha } \in C^\infty( \Lambda\times \bar\Omega; \R^{r\times m} )$ are
\begin{equation}
\label{lin1}
 a^{ij}_{\alpha }(\la,x) ={ \frac{\partial \cf_i}{\partial v_{j\alpha}}}(\la,x,w(x)),\ \text{and}\    b^{ij}_{\alpha }(\la,x, w(x)) ={ \frac{\partial \cg_i}{\partial v_{j\alpha}}}(\la,x,w(x)).
\end{equation}

By  \cite[Theorem 11.3 ]{[Pa]}, for each $\la \in  \Lambda$ and $w$ smooth, the Frechet derivative of the map $(f_\la, g_\la)$ 
at $w$  is the operator  \be\label{frechet} (L_{\la},B_{\la})\colon H^{k+s}(\Omega;\R^m ) \rightarrow  \hsbr \ee 
 induced on  Hardy-Sobolev spaces by  the differential operator \eqref{lin}.

A  differentiable map is {\it  semi-Fredholm}  if the Frechet derivative at any point is a linear semi-Fredholm operator.

\begin{proposition}\label{prop:ell1}  Let  $(\cf,\cg)$ be as above, with  $\cf(\la,x,0)=0,\,\cg(\la,x,0)=0.$  If,  for each $\la,$  the linearization   $(L_\la,B_\la)$  of  $(\cf_{\la},\cg_{\la})$  at $u\equiv 0$ is  elliptic and  $ s>n/2,$  then there exists an open  ball  $B=B(0,r) \subset \hkr $ such that the map 
\[  h=(f,g) \colon \Lambda \times B \rightarrow  H^{s}(\Omega;\R^m ) \times\hbbr \]
  induced by  $(\cf,\cg)$  is a continuously parametrized family of smooth semi-Fredholm maps. Moreover, if $\Lambda$  is a smooth manifold and $(\cf,\cg)$ is a smooth, then so is $h.$
  \end{proposition}

\proof  Since the estimate \eqref{schauder} holds, each $(L_\la,B_\la)$  is semi-Fredholm. On the other hand, the set of all semi-Fredholm operators is open. From this, by  compactness of $\Lambda,$ we can find a ball $B(0,r)$ such that 
$D(h_\la)(u)$  is semi-Fredholm for any $u\in B(0,r).$  This proves the first assertion. The second is clear. \qed\vskip5pt
\vskip 20pt
As a matter of fact,  under our assumptions, the map $h=(f,g)$ is a family of Fredholm maps.  This follows from the existence of  a rough parametrix of an elliptic boundary value problem \cite{[Ag-1],[W-R-L]}. While  the above proposition will be sufficient  for most of our needs, we will use the parametrix in order to prove that the set of bifurcation points of the family  $h$ arising in the proof of the theorem \ref{th:40} coincides with the set of bifurcation point of the elliptic system  \eqref{bvp1} in the sense of definition \ref{def:1}.   

We will be sketchy in what follows,  since the method is standard and we have only to notice  that the construction of a parametrix of an elliptic boundary value problem depends smoothly on parameters  (see \cite[Theorem 9.32]{[W-R-L]}, and  also \cite[Theorem 16.5]{[Ag]}, where boundary value problems for pseudo-differential operators with limited degree of smoothness are considered).

\begin{proposition}\label{prop:boots} Let the system \eqref{bvp1} verify the assumptions of the theorem \ref{th:40}, and let  $s> n/2.$ Then  the set $B$ of all bifurcation points of \eqref{bvp1}  in the sense of definition \ref{def:1} coincides with the set $Bif(h)$ of bifurcation points of the family  \[h \colon \R^q \times \hkr \rightarrow  \hsr \times\hbbr \]  defined by  \eqref{nbvp1}.\end{proposition}

\proof
Clearly $B\subset Bif(h).$  In order to prove the opposite inclusion we will use the standard elliptic bootstrap.
Keeping our previous notation, $ v_{i\alpha}$ (resp $ v'_{i\alpha}$)  will denote the components of a vector $ v\in  {\R}^{mk^*}$ (resp $ v' \in {\R}^{mk_{j}^*}$).

Since  $\cf(\la, x,0)=0,\,\cg(\la,x,0)=0,$  applying  \cite[Lemma 2.1] {[Mi-1]}  to each component of $\cf$ and to each $\cg_i$ we can write $(\cf,\cg)$ in the form:
\begin{equation}
\label{miln}
\begin{array}{l}
 {\cf}(\la,x,v) =\sum_{|\alpha |\leq m}
a_{\alpha }(\la, x,v) v_\alpha \\ 
\cg_i (\la,x,v') =\tau \sum_{|\alpha|\leq k_i}
 b^i_\alpha (\la,x, v')v'_{i\alpha };  1\leq i\leq r. 
\end{array}
\end{equation}
where $v_\alpha=(v_{1\alpha} \dots v_{m\alpha})^t,$ 

In order to simplify notations, we reparametrize  each  family  $\cB^i_{\la,v'}(x,D)$   by $v\in {\R}^{mk^*}$ using the  projectors $ \pi \colon {\R}^{mk^*}\to  {\R}^{mk_{i}^*}.$
In this way we obtain a family of boundary operators 
 $$\mathcal B_{\la,v}(x,D) = [\mathcal B_{\la,v}^1(x,D),\dots,\mathcal B_{\la,v}^r(x,D)]^t$$ parametrized by $\R^q\times {\R}^{mk^*}.$

Putting  $v=j_m(u))$ we have written
the map  $$(\cf,\cg) \colon \Lambda\times C^\infty(\bar \Omega;\R^m)\ra C^\infty(\bar \Omega;\R^m)\times C^\infty(\partial \Omega;\R^r) $$ 
in the form
\begin{equation}
\label{quasilin}
\begin{array}{l}
 {\cf}(\la,x,u, \dots, D^mu) =\cL_{\la,j_m(u)}(x,D) u
 \\
{\cg}(\la,x,u, \dots, D^{m}u) =\tau \cB_{\la, j_{m}(u)}(x,D) u.
\end{array}
\end{equation}
where $\cL,\cB$  linear differential operators depending on parameters $(\la,u).$

Now let us take  $v=0\in \R^{mk^*}$ and  observe  that,  by  \cite[Lemma 2.1] {[Mi-1]}, the pair \\ $(\cL_{\la,0}(x,D),\cB_{\la, 0}(x,D))$ coincides with the linearization  \eqref{lin} of the map $(\cf,\cg)$ at $u=0,$ which  is elliptic by hypothesis.  It follows from this that  for small enough  $\epsilon$ the restriction of the family  $\cH_{\la,v}(x,D) = (\cL_{\la,v}(x,D),\cB_{\la, v}(x,D)$ to  $\R^q\times B(0,\epsilon)\subset R^q\times \R^{mk*} $  is a family of elliptic boundary value problems. 

Let  us denote by 
$$H_{\la,v}= (L_{\la,v},B_{\la, v}) \colon\hkr \ra \hsbr$$   
 the operator induced by $\cH_{\la,v},$ on Hardy-Sobolev spaces. 

The construction of a rough parametrix  of an elliptic boundary value problem (see the proof  of  \cite[Theorem 9.32]{[W-R-L]}) uses  the inverse of the principal symbol, the canonical basis at points of the boundary and localization via smooth partitions of unity. Since  each of the above objects  behave well with respect to smooth variation of parameters, it follows  that any smooth family  of elliptic boundary value problems possesses a smooth  parametrix on a  neighborhood of a given point in the parameter space. 
 
Now, let  $\la_*\in Bif(h),$  and let $ \cP $ be  a left parametrix   of the family $\cH$ restricted to a neighborhood N of $(\la_*,0)$ in  $\R^q\times B(0,\epsilon).$

By definition, $ \cP $  is a family of operators  
 $$\cP_{\la,v}\colon C^\infty(\bar \Omega;\R^m)\times C^\infty(\partial \Omega;\R^r) \ra C^\infty(\bar \Omega;\R^m)$$ smoothly varying with $(\la,v)\in N$  which extends to  a smooth family of operators $$P_{\la,v}\colon \hsbr \ra \hkr$$  such that 
\begin{equation}\label{trix}
K_{\la,v}= P_{\la,v}H_{\la, v}-\Id_{\hkr} 
\end{equation}
  is a smooth family of bounded operators from $\hkr$ into ${H^{k+s+1}(\Omega;\R^m )}.$

Clearly the families $H_{\la,j^m(u)}, P_{\la,j^m(u)}$ and $K_{\la,j^m(u)}$ extend to families of bounded operators parametrized by a neighborhood $W$ of $(\la_*,0)$ in $\R^q\times\hkr.$

 Using \eqref{quasilin} we can rewrite  the restriction of $h$  to $W$  in the  form 
  \begin{equation}\label{fin}
h(\la,u)  = H_{\la, u} u
\end{equation} 
and therefore, by \eqref{trix}

  \begin{equation}\label{fine}
 P_{\la, u}h(\la,u)= u +  K_{\la, u}u
\end{equation} 
If  $(\la_n,u_n)\ra (\la_*,0)$ and $h(\la_n,u_n)=0,$
by \eqref{fine},  
$u_n = - K_{\la_n, u_n}u_n$ belongs to $H^{k+s+1}(\Omega;\R^m )$ and $u_n\ra 0$ in $H^{k+s+1}(\Omega;\R^m )$ as well.  Iterating this and using Sobolev embedding theorems
we obtain that $u_n\ra 0$ in  $C^{k}(\Omega;\R^m )$ for any $k,$  which proves that $\la_*$ belongs to $B.$
\qed\vskip5pt

 \vskip20pt
\subsection{C. Fedosov's  formula}\label{sec:5.3}
 Given a smooth manifold $M,$   $\he_c(M;\C)$ will  denote de Rham cohomology of complex valued  compactly supported forms of even degree. The  Chern-character is a natural transformation $ \ch \colon K_c(-) \ra \he_c(-;\C)$ preserving the  module structure over  the ring $K(-)$ and $\he(-;\C)$  respectively. 
If $\Lambda$ is a compact manifold,  the cohomological version of the Atiyah-Singer theorem for  families  $ \cS\colon \Lambda \rightarrow \pc$  states: 

\begin{equation}\label{HAS}
 \ch \Ind S = (-1)^n p_*(\ch [\sigma] )\, \hin  \he(\Lambda;\C). 
\end{equation}

Here $p_*$ is the push-forward homomorphism in de Rham cohomology called also {\it integration along the fiber.}   Integration along the fiber can be defined directly on differential forms.  Acting on compactly supported forms  on the total space of a smooth fiber bundle $\pi\colon E \to \Lambda $ with fiber $F,$  
the  {\it integration along the fiber}  $\oint_{F}$ is defined as follows:
 
 Let  us denote with $ \Omega_c^*(E)= \bigoplus_{i} \Omega_c^i(E)$  the smooth forms of mixed degree with compact support on $E.$   In  local coordinates $(\la_1,\dots,\la_q, x_1,\dots, x_n),$ where  the $\la$-s  are coordinates on the base  and  the $x$-s are coordinates on the fiber, we can write  a form  $ \theta \in \Omega_c^*(E)$    as $\theta = \theta' + \theta _{n},$ where  $\theta' $ contains all  terms  of degree less than $n$  in $dx_1, \dots, dx_{n}$ and 
 \[\theta _{n} =\sum_{i_1, \dots, i_r} f _{i_1 \dots i_r} (x,\la) dx_1\wedge, \dots, \wedge dx_{n} \wedge d\la_{i_1}\wedge \dots\wedge d\la_{i_r}.\] 
By definition,
\[ \oint_{F}\theta =   \oint_{F}\theta_{n}= \sum_{i_1 \dots i_r}[ \int_{F} f _{i_1 \dots i_r}(x,\la)\, dx_1\wedge \dots  \wedge dx_{n}]\, d\la_{i_1}\wedge \dots\wedge d\la_{i_r},\] 
where the integral inside the brackets  is the ordinary integral of a compactly supported form of maximal degree (see \cite{[Bo-Tu]}).

Using Chern-Weil theory of characteristic classes for smooth vector bundles over not necessarily compact manifolds   Fedosov obtained  an explicit expression for the smooth form representing the Chern character of the index bundle of  a  family of  pseudo-differential operators in $\pc$ in terms of its principal symbol.
 
 The following proposition is an immediate  consequence of  \cite[Corollary 6.5 ]{[Fe-1]}.

\begin{proposition}\label{pr:fed}
 If $ \mathcal{S}$ is a smooth family of  pseudo-differential operators in $\pc,$  then 
$\ch (\Ind S) = p_* \ch[\sigma] $ is the cohomology class of the form
$$ -\sum_{j=n}^\infty\displaystyle {\frac{(j-1)!} {(2\pi i)^j(2j-1)!}} \oint _{S^{2n-1}} tr(\sigma^{-1}d\sigma)^{2j-1},$$
where $S^{2n-1}=\partial B^{2n}$ is the boundary of a ball in $\R^{2n}$ such that the support of $\sigma$ is contained in  $\Lambda\times B^{2n}.$ 
\end{proposition}

\end{document}